\documentclass[11pt,reqno,a4paper]{article}
\usepackage{geometry}
\usepackage{amsmath,amsthm,amstext,amscd,amssymb,euscript,url}
\usepackage{mathrsfs}
\usepackage{calrsfs}
\usepackage[inline]{enumitem}
\usepackage{mathtools}
\usepackage{color}
\usepackage[hidelinks]{hyperref}
\usepackage[affil-it]{authblk}

\newcommand{\R}{\mathbb R}

\renewcommand{\phi}{\varphi}

\def\1{{\mathchoice {\rm 1\mskip-4mu l} {\rm 1\mskip-4mu l}
{\rm 1\mskip-4.5mu l} {\rm 1\mskip-5mu l}}}

\newtheorem{theorem}{{\small T}{\scriptsize HEOREM}}[section]
\newtheorem{corollary}[theorem]{{\bf{\small C}{\scriptsize OROLLARY}}}
\newtheorem{proposition}[theorem]{{\bf{\small P}{\scriptsize ROPOSITION}}}
\newtheorem{lemma}[theorem]{{\bf{\small L}{\scriptsize EMMA}}}
\newtheorem{remark}[theorem]{{\bf{\small R}{\scriptsize EMARK}}}
\newtheorem{definition}[theorem]{{\bf{\small D}{\scriptsize EFINITION}}}
\newtheorem{induction}[theorem]{{\bf{\small I}{\scriptsize NDUCTIVE HYPOTHESIS}}}
\newtheorem{assumption}[theorem]{{\bf{\small A}{\scriptsize SSUMPTION}}}

\renewcommand{\proofname}{\textup{\textbf{{\small P}{\scriptsize ROOF}}}}
\makeatletter
\renewenvironment{proof}[1][\proofname]{\par
  \pushQED{\qed}%
  \normalfont\noindent #1\@addpunct{.}\hspace{0.1cm}\ignorespaces
}{\popQED\par\@endpefalse}
\makeatother

\newcommand{\beq}{\begin{eqnarray}}
\newcommand{\eeq}{\end{eqnarray}}

\newcommand{\ba}{\begin{align*}}
\newcommand{\ea}{\end{align*}}

\newcommand{\be}{\begin{equation}}
\newcommand{\ee}{\end{equation}}

\newcommand{\bl}{\begin{lemma}}
\newcommand{\el}{\end{lemma}}

\newcommand{\br}{\begin{remark}}
\newcommand{\er}{\end{remark}}

\newcommand{\bt}{\begin{theorem}}
\newcommand{\et}{\end{theorem}}

\newcommand{\bd}{\begin{definition}}
\newcommand{\ed}{\end{definition}}

\newcommand{\bind}{\begin{induction}}
\newcommand{\eind}{\end{induction}}

\newcommand{\basu}{\begin{assumption}}
\newcommand{\esu}{\end{assumption}}

\newcommand{\bp}{\begin{proposition}}
\newcommand{\ep}{\end{proposition}}

\newcommand{\bc}{\begin{corollary}}
\newcommand{\ec}{\end{corollary}}

\newcommand{\bpr}{\begin{proof}}
\newcommand{\epr}{\end{proof}}

\newcommand{\bi}{\begin{itemize}}
\newcommand{\ei}{\end{itemize}}

\newcommand{\ben}{\begin{enumerate}}
\newcommand{\een}{\end{enumerate}}

\newmuskip\pFqmuskip

\newcommand\pFq[6][8]{%
	\begingroup 
	\pFqmuskip=#1mu\relax
	\mathcode`\,=\string"8000
	\begingroup\lccode`\~=`\,
	\lowercase{\endgroup\let~}\pFqcomma
	{}_{#2}F_{#3}{\left[\genfrac..{0pt}{}{#4}{#5};#6\right]}%
	\endgroup
}
\newcommand{\pFqcomma}{\mskip\pFqmuskip}

\title{{\bf Atomic limits of empirical measures of interacting particle systems with self-duality}}
\author[1]{Mario Ayala\thanks{Corresponding author: \href{mailto:mario.ayala@tum.de}{mario.ayala@tum.de}}}
\author[2]{Nicolas Dirr\thanks{\href{mailto:DirrNP@cardiff.ac.uk}{DirrNP@cardiff.ac.uk}}}
\author[1]{Johannes Zimmer\thanks{\href{mailto:jz@tum.de}{jz@tum.de}}}

\affil[1]{School of Computation, Information and Technology, Chair for Analysis and Modelling, Technische Universit\"at M\"unchen, Boltzmannstra\ss{}e 3, 85748 Garching, Germany}
\affil[2]{School of Mathematics, Cardiff University, Senghennydd Road, Cardiff CF24 4AG, United Kingdom}

\begin{document}

\maketitle

\begin{abstract}
We study the fixed-particle limit of empirical measures for a unified
self-dual family of conservative particle systems on \(\mathbb Z^d\): independent
random walkers, the symmetric exclusion process, and the symmetric inclusion
process. For a fixed number of particles \(k\), under diffusive space-time scaling parametrised by \(n\), we prove that as \(n\to\infty\) the normalised empirical measure converges in
\(D([0,T],\mathcal M_1(\mathbb R^d))\) to an atomic empirical measure carried by
\(k\) independent diffusions, with constant covariance determined by the
microscopic jump kernel.

Equivalently, the limit is the unique atomic solution of the Dean-Kawasaki
martingale problem with mass quantisation \(1/k\). In particular, exclusion and
inclusion have the same fixed-\(k\) atomic limit as independent walkers: the
interaction contributes only through microscopic near-collisions, whose
rescaled occupation time vanishes. The proof combines self-duality,
consistency-based reduction to two particles, and Green-function estimates for
finite perturbations of symmetric random walks.
\end{abstract}



\section{Introduction}

The large-scale behaviour of conservative interacting particle systems is
classically encoded in their empirical measure.  Depending on the normalisation,
however, this empirical measure may lead to rather different macroscopic
objects.  In the usual hydrodynamic regime one typically considers $O(n^d)$ particles in a
macroscopic box, normalises the empirical measure by $n^{-d}$, and rescales
diffusively.  Self-averaging then takes place: the empirical density converges
to an absolutely continuous measure whose density is the deterministic solution of a diffusion equation~\cite{demasi2006mathematical,kipnis1998scaling,ayala-valenzuela_hydrodynamic_2016}, and its $n^{-d/2}$
fluctuations converge to a Gaussian Ornstein-Uhlenbeck-type field
\cite{kipnis1998scaling,demasi2006mathematical,ayala2021higher}.  A different formal object is
suggested by fluctuating hydrodynamics, where the density field itself is
kept stochastic and carries multiplicative noise, leading formally to the
Dean-Kawasaki equation
\cite{dean1996langevin,kawasaki1994stochastic}.  These two viewpoints are often
close in spirit, but they are structurally different: in fluctuation theory the
drift is linearised around the hydrodynamic profile, by means of the
Boltzmann-Gibbs principle~\cite{rost1983hydrodynamik}, \cite{brox1984equilibrium},
\cite[Chapter~11]{kipnis1998scaling}, and the noise is evaluated at that
deterministic profile, hence the limit is Gaussian; in fluctuating hydrodynamics
both the drift and the noise amplitude are evaluated at the random density
itself.  For the self-dual family considered here, self-duality makes the drift
linear in the occupation variables, hence the two pictures differ only through
the noise; for non-linear gradient systems, such as zero-range processes, they
differ also in the drift.  For the self-dual family, a quantitative version of the
Boltzmann-Gibbs principle, with the order of every term made explicit, is
obtained by means of orthogonal polynomial
duality~\cite{ayala2018quantitative}, and the terms of higher order are the
subject of the second-order Boltzmann-Gibbs
principle~\cite{goncalves2017second}.
Quantitative central limit theorems for the density fluctuations of the simple
symmetric exclusion process, first at fixed times and then on path space with
higher-order lattice corrections, are obtained
in~\cite{gess2026quantitative,gess2026higher}.

Here we study a third scaling regime, complementary to both settings above: the fixed-$k$, or
atomic, regime. Here the number of particles is fixed, the empirical measure is
normalised by $1/k$ (and hence independent of the space-time rescaling parameter $n$),
and space and time are rescaled diffusively with $n$.
There is then no thermodynamic self-averaging, the total mass remains concentrated on
finitely many atoms, and the limiting object should remain stochastic.  The
natural candidate is the unregularised Dean-Kawasaki dynamics, understood not
as a density-valued SPDE but as a measure-valued martingale problem~\cite{konarovskyi2019dean}.

\subsection*{The particle systems and the main result}

We consider the standard self-dual family of symmetric conservative dynamics on
$\mathbb Z^d$: independent random walkers (IRW), the symmetric exclusion process
with maximum capacity $\alpha$ (SEP$(\alpha)$), and the symmetric inclusion process
(SIP$(\alpha)$). See Section~\ref{sec:microscopic-model} for the precise formulation.  These processes are encoded by a single interaction parameter
$\sigma\in\{0,-1,+1\}$ in~\eqref{eq:unified-k-generator}.  The common
self-duality structure~\cite{carinci2013duality,giardina2009duality} is used to
close the drift of the empirical measure exactly.

Self-duality enters here in a deliberately minimal way: all we use is that the
$k$-particle generator, applied to the occupation variable, reduces to the
one-particle generator, so that the drift of the empirical measure closes with
no reference to $\sigma$.  This is a small fragment of a considerably wider
theory, now recognised as a Lie-algebraic and quantum-group
phenomenon~\cite{giardina2026duality} and reaching well beyond the classical
duality functions and the symmetric rank-one setting, to orthogonal-polynomial
and dynamic dualities and to asymmetric, multi-species, and higher-spin models
\cite{franceschini2024orthogonal,wagenaar2026dynamic,kuan2026correlated}.

It is important to mention which part of this structure the present
argument actually uses, since the answer is not uniform across the proof.
The drift closure is purely algebraic: in the form we use it, it is the
statement that the one-particle duality function is carried by the dynamics
into itself, and one expects it to persist whenever a self-duality with a
one-particle duality function is available.  The other two ingredients are of a
different nature.  The reduction of the interaction term to a pair of particles
rests on consistency~\cite{carinci2021consistent}, and the vanishing of the
resulting collision term rests on reversibility and on a Nash inequality for
a symmetric finite-range walk.  It is this second group, and not the
duality, that confines the result below to the symmetric case.

Let $y^{(n),k}$ be the diffusively rescaled $k$-particle coordinate process and
set
\[
    \mu^{(n),k}_t
    =
    \frac1k\sum_{i=1}^k\delta_{y^{(n),k}_i(t)} .
\]
Our main result, Theorem~\ref{thm:empirical-to-dk}, states that for every fixed
$k$, for each of the three interactions, and for convergent deterministic
initial data,
\[
    \mu^{(n),k}\Longrightarrow \mu^k
    \qquad
    \text{in }D([0,T],\mathcal M_1(\mathbb R^d)),
\]
where
\[
    \mu^k_t
    =
    \frac1k\sum_{i=1}^k\delta_{Y_i(t)}
\]
and \(Y_1,\dots,Y_k\) are independent \(\mathbb R^d\)-valued diffusions with
constant-coefficient generator
\[
    A\phi
    =
    \frac{\alpha}{2}\,a:D^2\phi,
\]
where \(a=\sum_{r\in\mathbb Z^d}rr^\top p(r)\) is the constant covariance
matrix of the underlying jump kernel.
Equivalently, $\mu^k$ is the unique solution of the Dean-Kawasaki martingale
problem of Definition~\ref{def:DK-mp} with the corresponding atomic initial
condition.

Let us stress two points.  First, the limiting object is precisely the
atomic Dean-Kawasaki dynamics.  By the rigidity theorem of Konarovskyi,
Lehmann and von Renesse~\cite{konarovskyi2019dean}, the unregularised
Dean-Kawasaki martingale problem admits solutions only for atomic initial data
with quantised masses; in that case the solution is a system of independent
diffusions.  Hence our theorem gives a microscopic derivation of this rigid
object from finite-particle conservative systems.

Second, the limit is the same for IRW, SEP, and SIP: the interaction parameter
\(\sigma\) disappears.  The reason is that, at fixed \(k\), the interaction can
contribute only through microscopic near-collisions of labelled particles, and
the rescaled time spent on these diagonals vanishes.  Consequently the limiting
quadratic variation contains only the independent-walker contribution, namely
the mobility \(\alpha\rho\) evaluated on the atomic empirical measure.

This does not contradict hydrodynamic fluctuation theory.  In the hydrodynamic
regime the number of particles per macroscopic volume diverges, local
equilibrium averaging becomes relevant, and the mobilities are model-dependent:
\(\chi_{\mathrm{SEP}}(\rho)=\rho(\alpha-\rho)\) and
\(\chi_{\mathrm{SIP}}(\rho)=\rho(\alpha+\rho)\), against
\(\chi_{\mathrm{IRW}}(\rho)=\alpha\rho\).  These are density-level quantities:
they encode the variance of the microscopic current under local equilibrium
and therefore depend on the occupation density seen at a typical site.  In the
fixed-\(k\) atomic regime, by contrast, no macroscopic local equilibrium is
formed, the interaction part is confined to vanishing microscopic collision
times, and only the one-particle mobility survives in the limiting
Dean-Kawasaki bracket; see Section~\ref{subsec:loss-of-mobility}.

\subsection*{Main ideas of the proof}  For independent walkers the result is indeed elementary: the labelled coordinates are independent, each satisfies Donsker's invariance principle, and the convergence of the empirical measure follows from the continuity of the map \[     (y_1,\dots,y_k)\longmapsto \frac1k\sum_{i=1}^k\delta_{y_i}. \] The interacting cases are different.  For SEP and SIP the labelled law is not a product measure at any finite $n$: exclusion suppresses, and inclusion enhances, jumps onto occupied sites.  Hence one cannot simply apply Donsker's theorem to independent coordinates and then pass to empirical measures.  The non-trivial point is that these microscopic correlations are supported on near-collision events and leave no macroscopic trace in the fixed-$k$ scaling.  The proof has three ingredients.  First, self-duality closes the drift exactly: the one-particle duality function is the occupation variable, and the configuration generator acting on it reduces to the one-particle generator. Hence the drift of $\langle\mu^{(n),k},\phi\rangle$ is independent of $\sigma$ already at the microscopic level (Proposition~\ref{prop:drift_closure}).  Second, the interaction part of the carr\'e du champ is reduced to a two-particle quantity.  Since a fixed pair of particles in a $k$-particle SEP or SIP is not Markov in isolation, this reduction is not automatic.  We use the consistency intertwining of Carinci, Giardin\`a and Redig \cite{carinci2021consistent} to identify the relevant time-integrated $k$-particle near-collision functional with an average of genuine two-particle functionals (Lemma~\ref{lem:pairwise_reduction} and Proposition~\ref{prop:qv-two-particle-reduction}).  Third, the remaining two-particle term is a rescaled collision local time.  If \[     \Xi(t)=x^{(2)}_2(t)-x^{(2)}_1(t) \] is the two-particle difference process, then away from a finite set $\Xi$ is a symmetric finite-range random walk.  The interaction is therefore a finite local perturbation on $\mathbb Z^d$.  By means of a Nash inequality for the free walk, which is stable under finite perturbations, we obtain the on-diagonal decay of the transition probabilities of $\Xi$ (Proposition~\ref{prop:green-bounds-finite-perturb}), and from it, in Proposition~\ref{prop:vanish-exp-loc-time}, the following limit \[     \frac{1}{n^2}\int_0^{n^2T}\mathbb P(\Xi(u)=r)\,du     \longrightarrow 0 . \] This proves that the interaction contribution to the quadratic variation vanishes in $L^1$ (Corollary~\ref{cor:qv-interaction-vanishes}).  Together with tightness of the empirical measures (Theorem~\ref{thm:tightness-empirical-measures}), every subsequential limit solves the martingale problem of Definition~\ref{def:DK-mp} (Proposition~\ref{prop:identification-limit-points}), which by the uniqueness part of Theorem~\ref{thm:DK-uniqueness-rigidity} is well posed for atomic data; all subsequential limits therefore share one law and the whole sequence converges.  Consequently, although the interacting labelled processes are not independent microscopically, their fixed-$k$ empirical-measure limit has the same Dean-Kawasaki martingale problem as independent random walkers. 

\subsection*{Outline}

Section~\ref{sec:microscopic-model} introduces IRW, SEP and SIP,
the compatible configuration process, the basic self-duality relations, and the
fixed-$k$ empirical measure.  Section~\ref{sec:DK} recalls the formal
Dean-Kawasaki equation and formulates the associated martingale problem,
including its atomic solutions and rigidity.  The main convergence theorem is
then stated in Section~\ref{sec:main_result}, together with the interpretation
of the fixed-$k$ collapse.  The technical tools, the
consistency-based pair reduction and the Green-function estimates for finite
perturbations, are collected in Section~\ref{sec:technical-tools}.  The proof
of the convergence theorem is given in Section~\ref{sec:proof-convergence}.

\section{The microscopic model}
\label{sec:microscopic-model}
\subsection{The \texorpdfstring{$k$}{k}-particle coordinate process on \texorpdfstring{$\mathbb{Z}^d$}{Z\textasciicircum d}}
Fix $d\in\mathbb N$ and a symmetric jump kernel $p\colon \mathbb Z^d\to[0,1]$ satisfying
\[
p(r)=p(-r),\qquad \sum_{r\in\mathbb Z^d} p(r)=1.
\]
We assume that $p$ is \emph{irreducible}, i.e., the subgroup of $\mathbb Z^d$
generated by the set $\{r\in\mathbb Z^d:\ p(r)>0\}$ is all of $\mathbb Z^d$.
Moreover, we assume that $p$ is of \emph{finite range}, i.e., there exists
\(0<R<\infty\) such that \(p(r)=0\) whenever \(|r|>R\).

We write
\begin{equation}
\label{eq:def-second-moment-matrix}
    a:=\sum_{r\in\mathbb Z^d} rr^\top p(r)
\end{equation}
for the constant second-moment matrix of the jump kernel \(p\). Since \(p\) is
symmetric, $a$ is its covariance matrix.

\begin{remark}
\label{rem:a-positive-definite}
The matrix \(a\) is symmetric and positive definite. Indeed,
\(v^\top a\,v=\sum_{r}(r\cdot v)^2\,p(r)\ge0\) for \(v\in\mathbb R^d\), and if
\(v^\top a\,v=0\) then \(v\perp r\) for every \(r\) in the support of \(p\);
by irreducibility this support generates \(\mathbb Z^d\), hence spans
\(\mathbb R^d\), so \(v=0\). In particular the covariance matrix~\eqref{eq:def-second-moment-matrix} admits a symmetric
positive-definite square root \(a^{1/2}\), used below in the analysis of the
macroscopic limit.
\end{remark}

Let $k\in\mathbb N$ be fixed and write the coordinate vector
\[
x=(x_1,\dots,x_k)\in(\mathbb Z^d)^k, \qquad x_i\in\mathbb Z^d.
\]
For $i\in\{1,\dots,k\}$ and $r\in\mathbb Z^d$, define the updated configuration
\[
x^{i,i+r}:=(x_1,\dots,x_{i-1},x_i+r,x_{i+1},\dots,x_k).
\]
Introduce parameters $(\sigma,\alpha)\in \{0,1\}\times (0,\infty)\ \cup\ \{-1\}\times \mathbb{N}$
and define, for bounded $f\colon (\mathbb Z^d)^k\to\mathbb R$, the (permutation-invariant)
$k$-particle generator
\begin{equation}\label{eq:unified-k-generator}
(L^{(k)} f)(x)
:=\sum_{i=1}^k\ \sum_{r\in\mathbb Z^d} p(r)\,
\Big(\alpha+\sigma \sum_{\substack{j=1\\ j\neq i}}^k \mathbf 1_{\{x_j=x_i+r\}}\Big)\,
\big(f(x^{i,i+r})-f(x)\big).
\end{equation}
The corresponding coordinate process $x^{(k)}(\cdot)=(x_1(\cdot),\dots,x_k(\cdot))$
is the $(\mathbb Z^d)^k$-valued Markov process with generator $L^{(k)}$. When $\sigma=-1$ we implicitly restrict the process to initial data with
site occupancies bounded by $\alpha$ so that the jump rates in \eqref{eq:unified-k-generator} are nonnegative.
This constraint is preserved by the dynamics.
\medskip
The choices of \(\sigma\) correspond to the three standard models:
\begin{description}

\item[IRW.]
For \(\sigma=0\), particles are independent random walkers with jump rates
\(\alpha p(r)\) on \((\mathbb Z^d)^k\).

\item[SEP\((\alpha)\).]
For \(\sigma=-1\), with \(\alpha\in\mathbb N\), particles satisfy an exclusion
constraint with site capacity \(\alpha\). The state space is restricted to
\[
\sum_{j=1}^k \mathbf 1_{\{x_j=y\}} \le \alpha
\qquad \text{for all } y\in\mathbb Z^d .
\]
For \(\alpha=1\), this is the classical simple exclusion process
\cite{spitzer1970interaction}, equivalently restricted to
\[
E_k^{\mathrm{SEP}}
:=
\bigl\{x\in(\mathbb{Z}^d)^k:\ x_i\neq x_j\ \text{for all } i\neq j\bigr\}.
\]

\item[SIP\((\alpha)\).]
For \(\sigma=+1\), particles satisfy inclusion-type interaction: jumps to
occupied sites are enhanced by the arrival-site occupation.

\end{description}

\subsection{The compatible configuration process}
\label{subsec:compatible-config-process}

We also use the unlabelled configuration associated with the coordinate process.
Let
\[
\Omega
:=
\Big\{
\eta\in\mathbb N_0^{\mathbb Z^d}:
|\eta|:=\sum_{z\in\mathbb Z^d}\eta(z)<\infty
\Big\}.
\]
For \(x=(x_1,\dots,x_k)\in(\mathbb Z^d)^k\), set
\begin{equation}
\label{eq:def-Psi-main}
\Psi(x)
:=
\sum_{i=1}^k \delta_{x_i},
\qquad
\Psi(x)(z)=\sum_{i=1}^k \mathbf 1_{\{x_i=z\}} .
\end{equation}
The associated configuration process is
\[
    \eta^{(k)}(t):=\Psi\bigl(x^{(k)}(t)\bigr).
\]
It keeps track of occupation numbers and forgets particle labels. For \(u,v\in\mathbb Z^d\), \(u\neq v\), and \(\eta_u>0\), write
$\eta^{u,v}:=\eta-\delta_u+\delta_v $. The generator of the configuration process is
\begin{equation}
\label{eq:config-generator-main}
(\mathcal L f)(\eta)
=
\sum_{u,v\in\mathbb Z^d}
p(v-u)\,
\eta_u\,
\bigl(\alpha+\sigma\eta_v\bigr)
\bigl[f(\eta^{u,v})-f(\eta)\bigr],
\end{equation}
for bounded cylinder functions \(f\colon \Omega\to\mathbb R\), with the convention that
terms with \(u=v\) give no contribution. Equivalently,
\[
    (L^{(k)}(f\circ\Psi))(x)
    =
    (\mathcal L f)(\Psi(x)).
\]
Thus the choices \(\sigma=0,-1,1\) give, respectively, independent random
walkers, the exclusion process with capacity \(\alpha\), and the inclusion
process.  In the exclusion case \(\alpha\in\mathbb N\) and the state space is
restricted to
\[
\Omega_\alpha
:=
\{\eta\in\Omega:\eta_z\le \alpha \text{ for all }z\in\mathbb Z^d\}.
\]

\subsubsection{Self-duality}
\label{subsubsec:self-duality}

The configuration processes above are self-dual with respect to the standard duality functions for the IRW/SEP/SIP family~\cite{carinci2013duality,giardina2009duality}. We recall the form
needed in this paper.

For \(m,n\in\mathbb N_0\) let \((n)_m:=n(n-1)\cdots(n-m+1)\) denote the falling
factorial, with \((n)_0:=1\) and \((n)_m=0\) if \(m>n\), and let
\(d_\sigma(m,n):=(n)_m\,c_\sigma(m)\) be the one-site duality polynomials, with
normalisation \(c_\sigma(m)=\alpha^{-m}\) for \(\sigma=0\),
\(c_\sigma(m)=\Gamma(\alpha)/\Gamma(\alpha+m)\) for \(\sigma=1\), and
\(c_\sigma(m)=\Gamma(\alpha-m+1)/\Gamma(\alpha+1)\) for \(\sigma=-1\); in the
exclusion case the last expression is used only for \(0\le m\le\alpha\),
consistently with the capacity constraint. For finite configurations
\(\xi,\eta\in\Omega\) the self-duality function is
\begin{equation}
\label{eq:self-duality-function}
D(\xi,\eta)
:=
\prod_{x\in\mathbb Z^d}
d_\sigma(\xi_x,\eta_x),
\end{equation}
a finite product since \(\xi\) has finite support, which vanishes whenever
\(\xi_x>\eta_x\) at some site. Self-duality states that
\(\mathbb E_\eta[D(\xi,\eta(t))]=\mathbb E_\xi[D(\xi(t),\eta)]\) for all
\(t\ge0\) and all finite \(\xi,\eta\); equivalently, at generator level,
\((\mathcal L_\eta D)(\xi,\eta)=(\mathcal L_\xi D)(\xi,\eta)\), where
\(\mathcal L_\eta\) and \(\mathcal L_\xi\) act on the second and the first
argument respectively.

The one-particle case is especially important for the drift computation below.
For \(\xi=\delta_x\) and any \(\sigma\in\{-1,0,1\}\),
\[
D(\delta_x,\eta)
=
d_\sigma(1,\eta_x)
=
\frac{\eta_x}{\alpha}.
\]
Hence, up to the constant \(\alpha^{-1}\), the one-particle duality
function is simply the occupation variable; absorbing the constant into the
test function, we work with the unnormalised version
\(D_1(\delta_x,\eta):=\eta_x\). A direct computation gives
\begin{align}
(\mathcal L \eta_x)(\eta)
&=
\sum_{z,y\in\mathbb Z^d}
p(y-z)\eta_z(\alpha+\sigma\eta_y)
\bigl[
(\eta^{z,y})_x-\eta_x
\bigr]
\\
&=
\sum_{z\in\mathbb Z^d}
p(x-z)\eta_z(\alpha+\sigma\eta_x)
-
\sum_{y\in\mathbb Z^d}
p(y-x)\eta_x(\alpha+\sigma\eta_y).
\end{align}
The interaction terms cancel by the symmetry \(p(y-x)=p(x-y)\), leaving
\begin{equation}
\label{eq:occupation-duality-generator}
(\mathcal L \eta_x)(\eta)
=
\alpha\sum_{r\in\mathbb Z^d}p(r)
\bigl[\eta_{x+r}-\eta_x\bigr]
=
(L^{(1)}D_1(\delta_{\cdot},\eta))(x),
\end{equation}
where \(L^{(1)}\) is the one-particle coordinate generator
\((L^{(1)}\varphi)(x)=\alpha\sum_{r\in\mathbb Z^d}p(r)
[\varphi(x+r)-\varphi(x)]\).
This identity is the only self-duality result needed to close the drift of the
linear empirical observable in Proposition~\ref{prop:drift_closure}.

The full self-duality function \eqref{eq:self-duality-function} gives, for each
\(k\), a degree-\(k\) polynomial observable represented by the \(k\)-particle
coordinate process \eqref{eq:unified-k-generator}; we state it only for
context, since the quadratic variation below is closed through the consistency
reduction of Section~\ref{sec:technical-tools} rather than through this
function.

\subsection{The \texorpdfstring{$k$}{k}-particle empirical measure under diffusive rescaling}
\label{subsec:empirical-measure}
For $n\in\mathbb N$, define the diffusively rescaled coordinate process
\begin{equation}
\label{eq:def-rescY}
y^{(n),k}(t):=\frac{1}{n}x^{(k)}(n^2 t)\ \in\ (\mathbb R^d)^k,
\qquad\text{i.e.}\qquad
y^{(n),k}_i(t)=\frac{1}{n}x_i(n^2 t).
\end{equation}
Equivalently, $y^{(n),k}$ is a Markov process on $(\tfrac1n\mathbb Z^d)^k$ with generator
$\mathcal L^{(n),k}$ given by the diffusive time-change:
for any bounded function $f\colon (\tfrac1n\mathbb Z^d)^k \to \R$ and $y=x/n$,
\begin{equation}
\label{eq:diff-gen}
(\mathcal L^{(n),k} f)(y)=n^2\,(L^{(k)}F_n)(x),
\qquad\text{where}\qquad
F_n(x):=f(x/n).
\end{equation}

Let $\mathcal M_1(\mathbb R^d)$ denote the space of Borel probability
measures on $\mathbb R^d$, endowed with the topology of weak convergence,
metrised for instance by the Prokhorov metric. Since $\mathbb R^d$ is
Polish, $\mathcal M_1(\mathbb R^d)$ is Polish as well. For $T>0$, we
consider the Skorokhod space
\[
D([0,T],\mathcal M_1(\mathbb R^d)),
\]
equipped with the $J_1$ topology associated with this metric. For background
on weak convergence in metric spaces and on the Skorokhod topology in the
real-valued case, see Billingsley~\cite[Chaps.~1 and~3]{Billingsley1999};
for the general Polish-state-space setting used here, see
Ethier-Kurtz~\cite[Chap.~3, Secs.~5--8]{EthierKurtz1986}.

We define the empirical point measure induced by $y^{(n),k}$ as
\begin{equation}\label{eq:empirical-measure}
\mu^{(n),k}_t
:=\frac{1}{k}\sum_{i=1}^k \delta_{y^{(n),k}_i(t)}
\ \in \mathcal M_1(\mathbb R^d), \qquad t\ge 0.
\end{equation}
Since the coordinate process $y^{(n),k}$ is c\`adl\`ag, $\mu^{(n),k}$ is a c\`adl\`ag
$\mathcal M_1(\mathbb{R}^d)$-valued process and therefore defines a random
element of $D([0,T],\mathcal M_1(\mathbb{R}^d))$ for every $T>0$.

\begin{remark}\label{rem:kpart-vs-npart}
The normalisation in \eqref{eq:empirical-measure} places the empirical measure
$\mu^{(n),k}_t$ in a fundamentally different scaling regime from the one
used in classical hydrodynamic limits.

In the hydrodynamic setting (e.g., SEP or SIP at positive density) the
empirical density field is $\pi_t^{(n)}(\phi):=n^{-d}\sum_{x\in\mathbb Z^d}
\eta(x,n^2t)\phi(x/n)$, where $\eta(x,t)$ is the occupation number at site $x$;
the normalisation $n^{-d}$ compensates the diverging mass $\asymp\rho n^d$ of a
macroscopic box, and $\pi^{(n)}_t$ converges weakly to an absolutely continuous
measure $\pi_t(dx)=\rho(t,x)\,dx$ whose density solves
\begin{equation}
\label{eq:hydro-lim}
\partial_t \rho(t,x)
=
\frac{\alpha}{2}\,a:D^2\rho(t,x),
\end{equation}
with $a$ as in~\eqref{eq:def-second-moment-matrix};
see~\cite{kipnis1998scaling,demasi2006mathematical,ayala-valenzuela_hydrodynamic_2016}.

In contrast, the fixed-$k$ empirical measure \eqref{eq:empirical-measure} is
normalised by $1/k$, its total mass is constant in $n$, and the limiting object
remains purely atomic, $\mu_t^k=\frac1k\sum_{i=1}^k\delta_{y_i(t)}$. The two
constructions thus probe different scales, a thermodynamic, diffuse
regime versus a finite-particle, atomic one, and they share the same
linear macroscopic drift: in the atomic case
$\frac{d}{dt}\,\mathbb E[\langle\mu_t^k,\phi\rangle]
=\frac{\alpha}{2}\,\mathbb E[\langle\mu_t^k,a:D^2\phi\rangle]$ is a weak
formulation of \eqref{eq:hydro-lim}. The distinction lies not in the limiting
PDE but in the normalisation, and hence in the structure of the fluctuations.
\end{remark}

\section{The macroscopic model: the Dean-Kawasaki limit}
\label{sec:DK}

\subsection{Formal fluctuating hydrodynamics}

In this subsection we specialise to $d=1$ and to the independent case $\sigma=0$.
Let $p\colon \mathbb Z\to[0,1]$ be symmetric, irreducible and of finite range, with
scalar second moment $a$ as in~\eqref{eq:def-second-moment-matrix}.
For $n\in\mathbb N$ and $\phi\in C_c^\infty(\mathbb R)$ recall the diffusively rescaled
empirical density field
\begin{equation}\label{eq:pi-field-d1}
\pi_t^{(n)}(\phi)
:=
\frac{1}{n}\sum_{x\in\mathbb Z}\eta(x,n^2 t)\,
\phi\!\left(\frac{x}{n}\right),
\qquad t\ge 0.
\end{equation}
For each $\phi\in C_c^\infty(\mathbb R)$, Dynkin's formula yields the decomposition
\begin{equation}\label{eq:pi-martingale-decomp-d1}
\pi_t^{(n)}(\phi)
=
\pi_0^{(n)}(\phi)
+
\int_0^t \pi_s^{(n)}\!\left(\tfrac{\alpha a}{2}\phi''\right)\,ds
+
\int_0^t R_n\big(\eta(n^2 s),\phi\big)\,ds
+
M_t^{(n)}(\phi),
\end{equation}
where $M^{(n)}(\phi)$ is a c\`adl\`ag martingale and $R_n(\eta,\phi)$ is the
discretisation error arising from replacing the discrete generator by its
second-order Taylor approximation. Under the standing finite-range assumption on $p$,
one has the bound
\begin{equation}\label{eq:Rn-bound-formal}
\big|R_n(\eta,\phi)\big|
\le
\frac{C(\phi,p,\alpha)}{n}\,
\pi^{(n)}(\mathbf 1),
\qquad
\pi^{(n)}(\mathbf 1):=\frac{1}{n}\sum_{x\in\mathbb Z}\eta(x),
\end{equation}
so that $R_n=O(n^{-1})$ uniformly on time intervals where the macroscopic mass
$\pi^{(n)}(\mathbf 1)$ remains bounded.

\medskip
The predictable quadratic variation of $M^{(n)}(\phi)$ satisfies, formally,
\begin{equation}\label{eq:pi-bracket-leading}
\langle M^{(n)}(\phi)\rangle_t
=
\frac{\alpha a}{n}\int_0^t \pi_s^{(n)}(|\phi'|^2)\,ds
+ O\!\left(\frac{1}{n^2}\right).
\end{equation}
In particular, $M^{(n)}(\phi)$ is of size $n^{-1/2}$, which is the typical
fluctuation scale of the empirical density field in one dimension.

\medskip
Equations \eqref{eq:pi-martingale-decomp-d1}--\eqref{eq:pi-bracket-leading} suggest the weak
form of a fluctuating hydrodynamic equation: for each test function $\phi$,
\begin{equation}\label{eq:fh-weak-d1}
d\pi_t^{(n)}(\phi)
=
\pi_t^{(n)}\!\left(\tfrac{\alpha a}{2}\phi''\right)\,dt
+
dM_t^{(n)}(\phi),
\qquad
\mathbb E\!\left[(dM_t^{(n)}(\phi))^2\mid\mathcal F_t\right]
\approx
\frac{\alpha a}{n}\,\pi_t^{(n)}(|\phi'|^2)\,dt.
\end{equation}
Identifying the martingale term with a stochastic integral against a space-time
white noise $\dot W$ leads formally to the SPDE
\begin{equation}\label{eq:fh-spde-d1}
\partial_t \pi^{(n)}
=
\frac{\alpha a}{2}\,\partial_{xx}\pi^{(n)}
+
\frac{1}{\sqrt n}\,\partial_x\!\Big(\sqrt{\alpha a\,\pi^{(n)}}\,\dot W\Big),
\end{equation}
to be interpreted in the distributional sense, i.e., by testing against $\phi\in C_c^\infty(\mathbb R)$.

\subsection{The Dean-Kawasaki equation}
\label{subsec:DK-formal}

Equation \eqref{eq:fh-spde-d1} arises from the hydrodynamic scaling discussed in
Remark~\ref{rem:kpart-vs-npart}: the empirical density field is normalised by
$n^{-d}$, so that the total mass in a macroscopic region is of order one, while
its fluctuations are of order $n^{-d/2}$. In that regime the stochastic term
vanishes as $n\to\infty$, and one recovers the deterministic diffusion equation
\eqref{eq:hydro-lim}.

The situation is different in the fixed-$k$ regime of
the same remark. There, the total mass of the empirical measure
$\mu_t^{(n),k}$ is fixed and equal to $1$, independently of $n$, and no
self-averaging takes place. One therefore expects the macroscopic limit to remain
stochastic. Since the limiting object is an empirical measure of finitely many
diffusive particles, the natural candidate is a measure-valued equation whose drift
is given by the same second-order diffusion operator as in the hydrodynamic limit,
but whose quadratic variation reflects the finite-particle normalisation.

To motivate the correct form of the equation, consider first the limiting
atomic process $\mu_t^k=\frac1k\sum_{i=1}^k\delta_{Y_i(t)}$, where
$Y_1,\dots,Y_k$ are independent $\mathbb R^d$-valued diffusions with generator
$\mathcal A\phi=\frac{\alpha}{2}\,a:D^2\phi$. It\^o's formula shows
(Proposition~\ref{prop:DK-atomic-solutions} below) that
$\langle\mu_t^k,\phi\rangle$ has drift
$\frac{\alpha}{2}\langle\mu_t^k,a:D^2\phi\rangle$ and carries a continuous
martingale with quadratic variation
$\frac{\alpha}{k}\int_0^t\langle\mu_s^k,\nabla\phi\cdot a\nabla\phi\rangle\,ds$,
the factor $k^{-1}$ being the variance reduction from averaging $k$
independent one-particle martingales.

This calculation gives the formal Dean-Kawasaki dynamics in our setting.
If $\dot W=(\dot W_1,\dots,\dot W_d)$ denotes a standard
$\mathbb R^d$-valued space-time white noise and $a^{1/2}$ is the symmetric
square root of the covariance matrix $a$, then the corresponding formal SPDE is
\begin{equation}
\label{eq:DK-formal-corrected}
\partial_t \rho
=
\frac{\alpha}{2}\,a:D^2\rho
+
\frac{1}{\sqrt{k}}\,
\nabla\cdot\!\Big(\sqrt{\alpha\,\rho}\,a^{1/2}\dot W\Big);
\end{equation}
tested against $\phi\in C_c^\infty(\mathbb R^d)$ and formally integrated by
parts, it reproduces exactly the drift and quadratic variation above, i.e.,
those of the martingale problem in Definition~\ref{def:DK-mp} below.

Hence, in the fixed-$k$ regime, the Dean-Kawasaki equation is not obtained by
simply ``removing'' the hydrodynamic prefactor $n^{-d/2}$ from
\eqref{eq:fh-spde-d1}. Such a step would miss the normalisation of the empirical
measure and would produce the wrong noise strength. The correct coefficient is
instead dictated by the martingale structure of the atomic empirical measure
$\frac1k\sum_{i=1}^k\delta_{Y_i(t)}$.

Equation \eqref{eq:DK-formal-corrected} is still only formal, since for singular
measure-valued states the expression $\sqrt{\rho}\,\dot W$ has no direct meaning.
For this reason, throughout the remainder of this section we adopt the associated
martingale problem as the rigorous definition of Dean-Kawasaki dynamics. This is
also the natural formulation used in the well-posedness theory of~\cite{konarovskyi2019dean}, where it is proved that
solutions exist only in the integer regime and are necessarily atomic.

\subsection{The martingale formulation}
\label{subsec:DK-mp}

For $\phi\in C_c^\infty(\mathbb R^d)$, denote
\[
|\nabla \phi|_a^2 := \nabla \phi \cdot a \nabla \phi ,
\]
where $a$ is as in~\eqref{eq:def-second-moment-matrix}. 

\begin{definition}
\label{def:DK-mp}
Fix $k\in\mathbb N$ and an initial measure $\mu_0\in \mathcal M_1(\mathbb R^d)$.
A continuous $\mathcal M_1(\mathbb R^d)$-valued process
\[
\mu=(\mu_t)_{t\ge 0}
\]
is said to solve the Dean-Kawasaki martingale problem with parameter $k$
and initial condition $\mu(0)$ if $\mu_0=\mu(0)$ almost surely and, for every
$\phi\in C_c^\infty(\mathbb R^d)$,
\begin{equation}
\label{eq:DK-martingale}
M_t^\phi
:=
\langle \mu_t,\phi\rangle
-
\langle \mu_0,\phi\rangle
-
\frac{\alpha}{2}\int_0^t
\langle \mu_s,a:D^2\phi\rangle\,ds
\end{equation}
is a continuous martingale with quadratic variation
\begin{equation}
\label{eq:DK-bracket}
\langle M^\phi\rangle_t
=
\frac{\alpha}{k}\int_0^t
\langle \mu_s,|\nabla\phi|_a^2\rangle\,ds.
\end{equation}
\end{definition}

\subsection{Atomic solutions}
\label{subsec:DK-atomic}

The previous martingale problem is solved by empirical measures of independent
diffusions.

\begin{proposition}
\label{prop:DK-atomic-solutions}
Let $Y_1,\dots,Y_k$ be independent $\mathbb R^d$-valued diffusions with
generator
\[
\mathcal A\phi
=
\frac{\alpha}{2}\,a:D^2\phi,
\qquad \phi\in C_c^\infty(\mathbb R^d),
\]
and deterministic initial positions $y_1,\dots,y_k\in\mathbb R^d$.
Define
\begin{equation}
\label{eq:atomic-DK-solution}
\mu_t
:=
\frac{1}{k}\sum_{i=1}^k \delta_{Y_i(t)}.
\end{equation}
Then $(\mu_t)_{t\ge 0}$ is a solution to the martingale problem of
Definition~\ref{def:DK-mp} with initial condition
\[
\mu_0=\frac1k\sum_{i=1}^k \delta_{y_i}.
\]
\end{proposition}

\begin{proof}
For $\phi\in C_c^\infty(\mathbb R^d)$, It\^o's formula gives, for each $i$,
\[
\phi(Y_i(t))
=
\phi(Y_i(0))
+
\frac{\alpha}{2}\int_0^t a:D^2\phi(Y_i(s))\,ds
+
N_t^{i,\phi},
\]
where $N^{i,\phi}$ is a continuous martingale with
\[
\langle N^{i,\phi}\rangle_t
=
\alpha\int_0^t |\nabla\phi|_a^2(Y_i(s))\,ds.
\]
Averaging over $i=1,\dots,k$ yields
\[
\langle \mu_t,\phi\rangle
=
\langle \mu_0,\phi\rangle
+
\frac{\alpha}{2}\int_0^t
\langle \mu_s,a:D^2\phi\rangle\,ds
+
M_t^\phi,
\]
where
\[
M_t^\phi := \frac1k\sum_{i=1}^k N_t^{i,\phi}.
\]
Since the particles are independent, the martingales $N^{i,\phi}$ are orthogonal,
and therefore
\[
\langle M^\phi\rangle_t
=
\frac{1}{k^2}\sum_{i=1}^k \langle N^{i,\phi}\rangle_t
=
\frac{\alpha}{k}\int_0^t
\langle \mu_s,|\nabla\phi|_a^2\rangle\,ds.
\]
This is exactly \eqref{eq:DK-bracket}.
\end{proof}

\subsection{Uniqueness and rigidity of the Dean-Kawasaki martingale problem}
\label{subsec:DK-wellposedness}

For the convergence theorem below, the essential input is uniqueness in law for
the Dean-Kawasaki martingale problem when the initial condition is atomic with
masses \(1/k\). The same result also implies the stronger rigidity statement
that solutions exist only for such atomic initial data. This rigidity is a
special case of the theorem in~\cite{konarovskyi2019dean}. We remark that the rigidity holds for a wide class of equations of Dean-Kawasaki type~\cite{Muller2025a}. 

\begin{theorem}[\cite{konarovskyi2019dean}]
\label{thm:DK-uniqueness-rigidity}
Fix $k\in\mathbb N$ and $\mu_0\in\mathcal M_1(\mathbb R^d)$.
Then the Dean-Kawasaki martingale problem of
Definition~\ref{def:DK-mp} admits a solution if and only if
\begin{equation}
\label{eq:DK-atomic-initial}
\mu_0=\frac1k\sum_{i=1}^k \delta_{y_i}
\qquad\text{for some } y_1,\dots,y_k\in\mathbb R^d.
\end{equation}
In that case the solution is unique in law and is given by
\[
\mu_t=\frac1k\sum_{i=1}^k \delta_{Y_i(t)},
\]
where $Y_1,\dots,Y_k$ are independent diffusions with generator
$\frac{\alpha}{2}\,a:D^2$ and initial positions
$y_1,\dots,y_k$.
\end{theorem}

\begin{proof}
Existence for atomic initial data is Proposition~\ref{prop:DK-atomic-solutions}.
The converse and uniqueness in law follow
from~\cite[Theorem~2.2]{konarovskyi2019dean}, which is stated for a general
symmetric Markov diffusion operator. We apply it with the operator
\(\mathcal G:=a:D^2\) on \(\mathbb R^d\), whose carr\'e du champ is
\(\Gamma(\phi)=|\nabla\phi|_a^2\). By Remark~\ref{rem:a-positive-definite},
\(a\succ0\), so the semigroup of \(\mathcal G\) is Gaussian and commutes with
derivatives; hence it satisfies the two standing assumptions (2.2) and (2.3)
of~\cite{konarovskyi2019dean}, the second one with \(\rho=0\).

Let \((\mu_t)_{t\ge0}\) be a solution of the martingale problem of
Definition~\ref{def:DK-mp} with initial condition \(\mu_0\), and set
\(\nu_s:=\mu_{ks/\alpha}\). By \eqref{eq:DK-martingale} and
\eqref{eq:DK-bracket}, for every \(\phi\in C_c^\infty(\mathbb R^d)\) the
process \(s\mapsto M^\phi_{ks/\alpha}\) is a continuous martingale, and after
the change of variables in the time integrals it satisfies
\[
M^\phi_{ks/\alpha}
=
\langle\nu_s,\phi\rangle-\langle\nu_0,\phi\rangle
-\frac k2\int_0^s\langle\nu_u,\mathcal G\phi\rangle\,du,
\qquad
\langle M^\phi\rangle_{ks/\alpha}
=\int_0^s\langle\nu_u,\Gamma(\phi)\rangle\,du .
\]
Hence \((\nu_s)\) solves the martingale problem
of~\cite[Definition~2.1]{konarovskyi2019dean} with their parameter equal to
\(k\), tested against \(C_c^\infty(\mathbb R^d)\).

The class of test functions in~\cite{konarovskyi2019dean} is larger than
\(C_c^\infty(\mathbb R^d)\). Nevertheless, both identities extend to every
\(\phi\in C_b^2(\mathbb R^d)\), i.e., to
bounded functions with bounded continuous derivatives up to order two. Indeed,
if \(\chi\) is a smooth compactly supported cutoff equal to one near the
origin and \(\eta_\varepsilon\) is a standard mollifier, then
\(\phi_R:=\eta_{1/R}*(\chi(\cdot/R)\phi)\in C_c^\infty(\mathbb R^d)\), the
functions \(\phi_R\), \(\mathcal G\phi_R\) and \(\Gamma(\phi_R)\) are bounded
uniformly in \(R\), and they converge pointwise to \(\phi\),
\(\mathcal G\phi\) and \(\Gamma(\phi)\) as \(R\to\infty\). Since the \(\nu_u\)
are probability measures, the corresponding martingales are bounded uniformly
in \(R\) on every finite time interval. By dominated convergence we can
pass to the limit in the martingale identities for
\[
M^{\phi_R}_{ks/\alpha}
\quad\text{and}\quad
\bigl(M^{\phi_R}_{ks/\alpha}\bigr)^2
-\int_0^s\langle\nu_u,\Gamma(\phi_R)\rangle\,du,
\]
which gives both the martingale property and the stated quadratic variation
for \(\phi\).

This class is enough for our purposes. More precisely,
in~\cite{konarovskyi2019dean} the martingale problem enters only through the Laplace duality
of~\cite[Theorem~2.3]{konarovskyi2019dean}, whose proof tests it against the
functions \(-k\ln P_r e^{-f/k}\), \(r\in[0,t]\), where \((P_r)\) is the
semigroup of \(\tfrac k2\mathcal G\). For \(f\in C_b^2(\mathbb R^d)\) these
functions belong to \(C_b^2(\mathbb R^d)\), with bounds uniform in \(r\), since
\(P_re^{-f/k}\) is bounded below by a positive constant and \(P_r\) commutes
with derivatives.

By~\cite[Theorem~2.2]{konarovskyi2019dean}, a solution \((\nu_s)\) exists if
and only if \(\mu_0\) is of the form \eqref{eq:DK-atomic-initial}, in which
case it is unique in law and equals
\(\frac1k\sum_{i=1}^k\delta_{X_i(ks)}\), with \(X_1,\dots,X_k\) independent
diffusions with generator \(\tfrac12\mathcal G\) started at \(y_1,\dots,y_k\).
Notice that \(Y_i(t):=X_i(\alpha t)\) are independent diffusions with
generator \(\tfrac{\alpha}{2}\,a:D^2\). Hence, undoing the time change, we
obtain the stated representation
\[
\mu_t=\nu_{\alpha t/k}=\frac1k\sum_{i=1}^k\delta_{Y_i(t)} .\qedhere
\]
\end{proof}

\begin{remark}
\label{rem:DK-not-genuine-SPDE}
Theorem~\ref{thm:DK-uniqueness-rigidity} shows that, in the unregularised
setting considered here, the Dean-Kawasaki equation is not a genuinely
well-posed nonlinear SPDE on a large class of measure-valued initial data,
but a rigid reformulation of a finite system of independent diffusions: no
non-atomic solution exists. For the present paper, this rigidity identifies precisely the correct limiting object
for the fixed-$k$ empirical measures $\mu_t^{(n),k}$, and our convergence
theorem should be read as a microscopic derivation of the unique admissible,
namely atomic, Dean-Kawasaki dynamics in this regime.
\end{remark}

\section{Convergence of Empirical Measures to Dean-Kawasaki}
\label{sec:main_result}

\subsection{Main result: convergence to the atomic Dean-Kawasaki dynamics}

We consider the empirical measure process $\mu^{(n),k}$ of \eqref{eq:empirical-measure} as a random element of the space
$D([0,T],\mathcal M_1(\mathbb R^d))$, with the topology of
Section~\ref{subsec:empirical-measure}. Its candidate limit is the continuous
 $\mathcal M_1(\mathbb{R}^d)$-valued process
$\mu^k_t=\frac{1}{k}\sum_{i=1}^k \delta_{Y_i(t)}$, where $Y_1,\dots,Y_k$ are
independent $\mathbb{R}^d$-valued diffusions with generator
$A\varphi=\frac{\alpha}{2}a:D^2\varphi$; by
Proposition~\ref{prop:DK-atomic-solutions}, it solves the Dean-Kawasaki
martingale problem of Definition~\ref{def:DK-mp}.

The following result identifies the fixed-$k$ scaling limit of the microscopic
empirical measure process. In this regime the limit is not a diffuse
density-valued object, but the atomic Dean-Kawasaki dynamics generated by
$k$ independent limiting particles.

\begin{theorem}
\label{thm:empirical-to-dk}
Fix $k\in\mathbb N$ and let $y^{(n),k}$ be the diffusively rescaled
coordinate process defined in \eqref{eq:def-rescY}. Let
\[
\mu^{(n),k}_t
=
\frac{1}{k}\sum_{i=1}^k \delta_{y^{(n),k}_i(t)}
\]
denote the associated empirical measure process. Assume that
\[
y^{(n),k}(0)\to y^k=(y_1,\dots,y_k)
\quad
\text{in }(\mathbb R^d)^k .
\]
For $\sigma=-1$ we also assume that every microscopic initial configuration is
admissible, i.e., it has at most $\alpha$ particles per site; the limiting
positions $y_1,\dots,y_k$ may coincide.
Then, for every $T>0$,
\[
\mu^{(n),k} \Rightarrow \mu^k
\quad
\text{in }
D([0,T],\mathcal M_1(\mathbb R^d)),
\]
where
\[
\mu^k_t=\frac{1}{k}\sum_{i=1}^k\delta_{Y_i(t)}
\]
and $(Y_i)_{i=1}^k$ are independent diffusions with generator
\[
A\varphi
=
\frac{\alpha}{2}a:D^2\varphi.
\]
Equivalently, the limiting process $\mu^k$ is the unique solution of the
Dean-Kawasaki martingale problem of Definition~\ref{def:DK-mp} with initial
condition
\[
\mu^k_0
=
\frac{1}{k}\sum_{i=1}^k\delta_{y_i}.
\]
\end{theorem}

\subsection{Interpretation: loss of the interaction mobility in the fixed-\texorpdfstring{$k$}{k} limit}
\label{subsec:loss-of-mobility}
In the proof of Theorem~\ref{thm:empirical-to-dk} we can see how the parameter $\sigma$ disappears in the
fixed-$k$ limit: the interaction
contribution to the microscopic carr\'e du champ, the term
\(R_{\mathrm{int}}^{(n),k}\) of \eqref{eq:Rint-def} below, is supported on
microscopic near-collisions \(y_j=y_i+r/n\), and its time integral vanishes in
expectation (Corollary~\ref{cor:qv-interaction-vanishes}). 

What controls the survival of $\sigma$ is not $k$ itself but the density
$\rho=k/n^d$ of particles in a macroscopic box, since $\sigma$ can reach the
limit only through the collision term. This is the regime studied
in~\cite{ayala2025collisions}, where collision counts and unidirectional fluxes
of the exclusion process are retained as observables in their own right under a
joint scaling of the lattice spacing and the particle number. There the
collision observable has a genuinely regime-dependent limit: in the sparse
regime collisions become rare and the exclusion dynamics approaches that of
independent walkers, whereas at hydrodynamic density they persist and are
responsible for the nonlinear mobility. Theorem~\ref{thm:empirical-to-dk} is
the extreme case of that picture, with $k$ fixed.

\section{Technical Tools}
\label{sec:technical-tools}
\subsection{Consistency and reduction to two particles}
We use consistency of the coordinate processes to reduce averaged
\(k\)-particle observables to two-particle observables.
Let \(V=\mathbb Z^d\) be the underlying site space, and denote by
\(P_t^{(k)}\) the transition semigroup of the \(k\)-particle coordinate
process. For \(k\ge2\), define the random-removal operator by
\[
(\Pi^{(k)}g)(x_1,\ldots,x_k)
:=
\frac1k\sum_{i=1}^k
g(x_1,\ldots,x_{i-1},x_{i+1},\ldots,x_k).
\]
Thus, \(\Pi^{(k)}g\) evaluates \(g\) after deleting one coordinate
chosen uniformly at random.

The relevant consistency property is the random-removal intertwining
of Carinci--Giardin\`a--Redig~\cite{carinci2021consistent}.
In our notation, it reads
\begin{equation}
\label{eq:semigroup-intertwining-from-consistency}
P_t^{(k)}\Pi^{(k)}g
=
\Pi^{(k)}P_t^{(k-1)}g,
\qquad t\ge0,
\end{equation}
for every bounded permutation-invariant function
\(g\colon V^{k-1}\to\mathbb R\).
This identity says that evolving the \(k\)-particle system and then
removing a uniformly chosen particle gives the same expectation
as removing a uniformly chosen particle first and then evolving
the remaining \(k-1\) particles.

The following proposition shows that this hypothesis is met by all three
models treated in this paper.

\begin{proposition}
\label{prop:models-consistent}
For each admissible parameter pair \((\sigma,\alpha)\), the configuration
process \eqref{eq:config-generator-main}, compatible with the \(k\)-particle
coordinate processes generated by \eqref{eq:unified-k-generator}, is consistent
in the sense of Definition~3.1 of~\cite{carinci2021consistent}. Hence the semigroup intertwining
\eqref{eq:semigroup-intertwining-from-consistency} holds for every bounded
permutation-invariant \(g\colon V^{k-1}\to\mathbb R\).
\end{proposition}

\begin{proof}
In the compatible configuration process \eqref{eq:config-generator-main}, with \(\sigma \in \{0,-1,1\} \), a
particle jumps from \(u\) to \(v\) at rate
\(\eta_u\,\bigl(\sigma p(v-u)\,\eta_v+\alpha p(v-u)\bigr)\), which is of the form
(3.6) in~\cite[Theorem~3.3]{carinci2021consistent} with
\(\theta(\{u,v\})=\sigma p(v-u)\) and \(\alpha(u,v)=\alpha p(v-u)\). Notice that
this theorem is stated for a countable vertex set, hence it applies to
\(\mathbb Z^d\), and we obtain that the generator of the configuration process
commutes with the annihilation operator, i.e., the configuration process is
consistent.
By~\cite[Theorem~2.7]{carinci2021consistent} this is equivalent to the
generator-level intertwining \(L^{(k)}\Pi^{(k)}=\Pi^{(k)}L^{(k-1)}\) on
permutation-invariant functions, which is equation (2.14) there.

It remains to pass from generators to semigroups. For fixed \(k\) the total
jump rate of \eqref{eq:unified-k-generator} is bounded: each particle jumps at
rate at most \((\alpha+k-1)\sum_r p(r)=\alpha+k-1\), so
\(\|L^{(k)}\|_\infty\le 2k(\alpha+k-1)<\infty\) on bounded functions, and
\(P_t^{(k)}=e^{tL^{(k)}}\) is given by a norm-convergent exponential series. Since
\(L^{(k-1)}\) preserves permutation-invariance, iterating
\(L^{(k)}\Pi^{(k)}=\Pi^{(k)}L^{(k-1)}\) gives
\((L^{(k)})^n\Pi^{(k)}=\Pi^{(k)}(L^{(k-1)})^n\) for all \(n\ge0\), whence
\[
P_t^{(k)}\Pi^{(k)}g
=\sum_{n\ge0}\frac{t^n}{n!}\,(L^{(k)})^n\Pi^{(k)}g
=\sum_{n\ge0}\frac{t^n}{n!}\,\Pi^{(k)}(L^{(k-1)})^ng
=\Pi^{(k)}P_t^{(k-1)}g,
\]
which is \eqref{eq:semigroup-intertwining-from-consistency}.
\end{proof}

For bounded \(h\colon V^2\to\mathbb R\) and \(k\ge 2\), define the ordered-pair
symmetrisation
\begin{equation}
\operatorname{OSym}_k(h)(x_1,\dots,x_k)
\,:=\,
\frac{1}{k(k-1)}
\sum_{\substack{1\le i,j\le k\\ i\neq j}}
 h(x_i,x_j).
\end{equation}

\noindent
Then we have the following.

\begin{lemma}
\label{lem:pairwise_reduction}
Assume that the family \(\big(X^{(k)}\big)_{k\ge 1}\) is permutation-invariant and satisfies
the consistency intertwining \eqref{eq:semigroup-intertwining-from-consistency}.
Then, for every \(k\ge 2\), every \(t\ge 0\), and every bounded
\(h\colon V^2\to\mathbb R\),
\begin{equation}
P_t^{(k)}\,\operatorname{OSym}_k(h)
=
\operatorname{OSym}_k\!\big(P_t^{(2)}h\big).
\label{eq:OSym_semigroup_identity}
\end{equation}
Equivalently, for every deterministic initial condition
\(x=(x_1,\dots,x_k)\in V^k\),
\begin{equation}
\mathbb E_x^{(k)}
\left[
\operatorname{OSym}_k(h)(X^{(k)}(t))
\right]
=
\frac{1}{k(k-1)}
\sum_{\substack{1\le i,j\le k\\ i\neq j}}
\mathbb E_{(x_i,x_j)}^{(2)}
\left[
h(X^{(2)}(t))
\right].
\label{eq:pair_reduction_deterministic_initial}
\end{equation}
\end{lemma}

\begin{proof}
For $k=2$, write $h^\leftarrow(u,v):=h(v,u)$, so that
\[
\operatorname{OSym}_2(h)
=
\frac12\big(h+h^\leftarrow\big).
\]
By permutation invariance of the two-particle dynamics, $P_t^{(2)}$ commutes
with the exchange of the two labels, i.e.,
$P_t^{(2)}(h^\leftarrow)=(P_t^{(2)}h)^\leftarrow$. Hence we obtain
\[
P_t^{(2)}\operatorname{OSym}_2(h)
=
\frac12\big(P_t^{(2)}h+(P_t^{(2)}h)^\leftarrow\big)
=
\operatorname{OSym}_2\big(P_t^{(2)}h\big).
\]
For $k\ge3$, observe first that $\operatorname{OSym}_{k-1}(h)$ is
permutation-invariant. Moreover,
\begin{equation}\label{eq:Sym_perm_inv}
 \Pi^{(k)}\operatorname{OSym}_{k-1}(h)
=
\operatorname{OSym}_k(h).   
\end{equation}
Indeed, each ordered pair $(i,j)$ with $i\neq j$ appears exactly for the
$k-2$ choices of the removed label distinct from $i$ and $j$, and the
normalisations give
\[
\frac1k \cdot \frac{k-2}{(k-1)(k-2)}
=
\frac1{k(k-1)}.
\]
Using~\eqref{eq:Sym_perm_inv}, the intertwining relation~\eqref{eq:semigroup-intertwining-from-consistency}, we find
\[
P_t^{(k)}\operatorname{OSym}_k(h)
=
P_t^{(k)}\Pi^{(k)}\operatorname{OSym}_{k-1}(h)
=
\Pi^{(k)}P_t^{(k-1)}\operatorname{OSym}_{k-1}(h).
\]
Iterating this identity down to two particles yields
\[
P_t^{(k)}\operatorname{OSym}_k(h)
=
\Pi^{(k)}\Pi^{(k-1)}\cdots \Pi^{(3)}
P_t^{(2)}\operatorname{OSym}_2(h).
\]
Applying the elementary identity
\[
\Pi^{(k)}\Pi^{(k-1)}\cdots \Pi^{(3)}\operatorname{OSym}_2(f)
=
\operatorname{OSym}_k(f)
\]
with $f=P_t^{(2)}h$ gives \eqref{eq:OSym_semigroup_identity}. Evaluating at
$x=(x_1,\dots,x_k)$ gives \eqref{eq:pair_reduction_deterministic_initial}.
\end{proof}

\begin{remark}
Lemma~\ref{lem:pairwise_reduction} is an averaged pair-reduction
identity: in \eqref{eq:pair_reduction_deterministic_initial}, each
two-particle expectation on the right-hand side is started from the ordered
pair \((x_i,x_j)\) extracted from \(x\). It should not be read as saying that
a fixed pair of coordinates of \(X^{(k)}\) evolves as an autonomous
two-particle process: the expectations on the right-hand side are computed with
the two-particle semigroup, started from the different initial pairs.
\end{remark}

\begin{corollary}
\label{cor:time_integrated_pair_reduction}
Assume the hypotheses of Lemma~\ref{lem:pairwise_reduction}. Let
$h_t\colon V^2\to\mathbb R$, $t\in[0,T]$, be a bounded measurable family of
functions. Then, for every deterministic initial condition
$x=(x_1,\dots,x_k)\in V^k$,
\[
\mathbb E_x^{(k)}
\left[
\int_0^T
\operatorname{OSym}_k(h_t)(X^{(k)}(t))\,dt
\right]
=
\frac{1}{k(k-1)}
\sum_{\substack{1\le i,j\le k\\ i\neq j}}
\int_0^T
\mathbb E_{(x_i,x_j)}^{(2)}
\left[
h_t(X^{(2)}(t))
\right]dt .
\]
\end{corollary}

\begin{proof}
Apply Lemma~\ref{lem:pairwise_reduction} to \(h_t\) for each fixed
\(t\in[0,T]\), integrate in \(t\), and use Fubini's theorem.
\end{proof}

\subsection{Truncated Green bounds for finite perturbations}

We next introduce the estimate needed to control two-particle near-collisions. It
applies to reversible Markov chains on \(\mathbb Z^d\) which, outside a finite
set, coincide with a symmetric irreducible finite-range random walk. The result
bounds the truncated Green function uniformly in the starting point, with the
usual dimension-dependent growth: order \(\sqrt t\) in \(d=1\), logarithmic
growth in \(d=2\), and a uniform bound in \(d\ge3\). By means of a comparison
of Dirichlet forms, we transfer a Nash inequality for the free walk to the
perturbed chain, and from it we obtain the decay of the transition
probabilities.

For \(S\subseteq\mathbb Z^d\) and \(f\colon S\to\mathbb R\) we write \(\|f\|_1\) and
\(\|f\|_2\) for the norms of \(f\) in \(\ell^1(S)\) and \(\ell^2(S)\) with
respect to the counting measure, and \(\|f\|_{1,m}\) and \(\|f\|_{2,m}\) for the
corresponding norms with respect to a measure \(m\) on \(S\).

\begin{lemma}
\label{lem:nash-free-walk}
Let \(d\ge1\), let \(q_0\colon \mathbb Z^d\to[0,\infty)\) be a symmetric finite-range
jump kernel whose support generates \(\mathbb Z^d\), and let \(S\) be either
\(\mathbb Z^d\) or \(\mathbb Z^d\setminus\{v\}\) for some \(v\in\mathbb Z^d\).
For finitely supported \(f\colon S\to\mathbb R\) define
\[
\mathcal E_0^S(f,f)
:=
\frac12\sum_{z,w\in S}q_0(w-z)\bigl(f(w)-f(z)\bigr)^2 .
\]
Then there exists a constant \(C_N<\infty\), depending only on \(d\) and
\(q_0\), such that, for every finitely supported \(f\colon S\to\mathbb R\),
\begin{equation}
\label{eq:nash-free-walk}
\|f\|_2^{2+4/d}
\le
C_N\,\mathcal E_0^S(f,f)\,\|f\|_1^{4/d}.
\end{equation}
\end{lemma}

\begin{proof}
Since the support of \(q_0\) generates \(\mathbb Z^d\), it spans
\(\mathbb R^d\), hence it contains \(d\) linearly independent vectors
\(b_1,\dots,b_d\). Let \(H:=\mathbb Zb_1+\dots+\mathbb Zb_d\), which is a
subgroup of finite index in \(\mathbb Z^d\), and define
\[
D^S(f)
:=
\sum_{j=1}^d\ \sum_{\substack{z\in S\\ z+b_j\in S}}
\bigl(f(z+b_j)-f(z)\bigr)^2 .
\]
Keeping in \(\mathcal E_0^S\) only the terms with
\(w-z\in\{\pm b_1,\dots,\pm b_d\}\), we obtain the following lower bound
\begin{equation}
\label{eq:energy-retained-edges}
\mathcal E_0^S(f,f)\ge c_0\,D^S(f),
\qquad
c_0:=\min_{1\le j\le d}q_0(b_j)>0 .
\end{equation}
Notice that the edges \(\{z,z+b_j\}\) retained in \(D^S\) only connect points
of the same coset of \(H\), and that, for every \(x\in\mathbb Z^d\), the map
\((n_1,\dots,n_d)\mapsto x+n_1b_1+\dots+n_db_d\) is a graph isomorphism between
the nearest-neighbour lattice \(\mathbb Z^d\) and the coset \(x+H\) endowed
with these edges. If \(S=\mathbb Z^d\setminus\{v\}\), the coset \(v+H\) loses
the vertex \(v\) and becomes a copy of the punctured lattice
\(\mathbb Z^d\setminus\{0\}\), and all the other cosets are unchanged.

Hence it is enough to combine the finitely many cosets, once we know that, for
\(G=\mathbb Z^d\) and for \(G=\mathbb Z^d\setminus\{0\}\), every finitely
supported \(g\colon G\to\mathbb R\) satisfies the nearest-neighbour Nash
inequality
\begin{equation}
\label{eq:nash-nn}
\|g\|_2^{2+4/d}
\le
C_d\,D_G(g)\,\|g\|_1^{4/d},
\qquad
D_G(g):=\sum_{\substack{\{x,y\}\subset G\\ |x-y|_1=1}}\bigl(g(x)-g(y)\bigr)^2 .
\end{equation}
We now prove the Nash inequality. We consider \(d=1\), so that \(G\) is either \(\mathbb Z\) or the union of two
half-lines. In both cases every \(x\in G\) is the starting point of an infinite
nearest-neighbour path \(x=x_0,x_1,x_2,\dots\) in \(G\) without repeated
vertices. Since \(g\) is finitely supported, by the Cauchy-Schwarz inequality
we have
\[
g(x)^2
=
\sum_{i\ge0}\bigl(g(x_i)-g(x_{i+1})\bigr)\bigl(g(x_i)+g(x_{i+1})\bigr)
\le
2\,D_G(g)^{1/2}\,\|g\|_2 .
\]
Together with \(\|g\|_2^2\le\|g\|_\infty\|g\|_1\) we obtain
\(\|g\|_2^4\le2\,D_G(g)^{1/2}\|g\|_2\|g\|_1^2\), i.e.,
\(\|g\|_2^6\le4\,D_G(g)\|g\|_1^4\), which is \eqref{eq:nash-nn} for \(d=1\).

Now the case \(d\ge2\). For \(G=\mathbb Z^d\) the inequality \eqref{eq:nash-nn} is
classical: the transition probabilities of the continuous-time simple random
walk on \(\mathbb Z^d\) are bounded by \(c\,t^{-d/2}\) for all \(t>0\), and
\eqref{eq:nash-nn} follows from the converse implication
in~\cite[Theorem~(2.1)]{carlen1987upper}; see the proof
of~\cite[Proposition~3.1]{bass2008symmetric}.

For \(G=\mathbb Z^d\setminus\{0\}\) and \(g\colon G\to\mathbb R\) finitely supported, we
extend \(g\) to \(\hat g\colon \mathbb Z^d\to\mathbb R\) by setting
\(\hat g(0):=g(e_1)\), so that \(\|\hat g\|_2\ge\|g\|_2\) and
\(\|\hat g\|_1\le2\|g\|_1\). The two energies differ only on the edges
incident to the origin, i.e.,
\[
D_{\mathbb Z^d}(\hat g)
=
D_G(g)+\sum_{|u|_1=1}\bigl(g(u)-g(e_1)\bigr)^2 .
\]
Every neighbour \(u\) of the origin is joined to \(e_1\) by a nearest-neighbour
path in \(G\) with at most four edges, namely \(e_1,\,e_1+u,\,u\) for
\(u=\pm e_j\) with \(j\ge2\), and
\(e_1,\,e_1+e_2,\,e_2,\,-e_1+e_2,\,-e_1\) for \(u=-e_1\). Hence
\(\bigl(g(u)-g(e_1)\bigr)^2\le4\,D_G(g)\) by the Cauchy-Schwarz inequality, and
therefore \(D_{\mathbb Z^d}(\hat g)\le(1+8d)\,D_G(g)\). Applying
\eqref{eq:nash-nn} on \(\mathbb Z^d\) to \(\hat g\), we obtain
\eqref{eq:nash-nn} on \(\mathbb Z^d\setminus\{0\}\) with a larger constant.

We just need to combine the cosets as follows. Let \(S_1,\dots,S_N\) be the intersections
of \(S\) with the cosets of \(H\), let \(f_i\) be the restriction of \(f\) to
\(S_i\), and let \(D_i(f)\) be the part of \(D^S(f)\) coming from the edges
inside \(S_i\), so that \(D^S(f)=\sum_iD_i(f)\). By \eqref{eq:nash-nn} and the
isomorphism above,
\(\|f_i\|_2^2\le C\,D_i(f)^{d/(d+2)}\|f_i\|_1^{4/(d+2)}\) for every \(i\), with
a constant \(C\) depending only on \(d\). Summing over \(i\), and using
H\"older's inequality with exponents \((d+2)/d\) and \((d+2)/2\) together with
the bound \(\sum_i\|f_i\|_1^2\le\|f\|_1^2\), we obtain
\[
\|f\|_2^2
\le
C\Bigl(\sum_{i=1}^ND_i(f)\Bigr)^{\frac{d}{d+2}}
\Bigl(\sum_{i=1}^N\|f_i\|_1^2\Bigr)^{\frac{2}{d+2}}
\le
C\,D^S(f)^{\frac{d}{d+2}}\,\|f\|_1^{\frac{4}{d+2}} .
\]
Raising this inequality to the power \((d+2)/d\) and using
\eqref{eq:energy-retained-edges}, we obtain \eqref{eq:nash-free-walk}.
\end{proof}

\begin{proposition}
\label{prop:green-bounds-finite-perturb}
Let \(d\ge 1\), and let \(Z=(Z(t))_{t\ge0}\) be a continuous-time
reversible Markov chain on \(\mathbb Z^d\) with generator
\[
(L_Z f)(z)=\sum_{w\in\mathbb Z^d} q(z,w)\bigl(f(w)-f(z)\bigr),
\qquad z\in\mathbb Z^d.
\]
We take \(q(z,z)=0\) and assume moreover that:

\begin{enumerate}[label=(\roman*)]
\item\label{hyp:green-bounded-rates} \(Z\) has bounded total jump rates,
\[
\sup_{z\in\mathbb Z^d}\sum_{w\in\mathbb Z^d} q(z,w)<\infty,
\]
and finite jump range;

\item\label{hyp:green-finite-perturbation} there exist a finite set \(F\subset\mathbb Z^d\) and a symmetric
finite-range jump kernel \(q_0\colon \mathbb Z^d\to[0,\infty)\) whose support
generates \(\mathbb Z^d\) such that, for every \(z\notin F\),
\[
q(z,z+r)=q_0(r),\qquad r\in \mathrm{supp}(q_0)\setminus\{0\},
\]
and \(q(z,w)=0\) otherwise;

\item\label{hyp:green-support-reversible} the rates \(q\) have the same support as the free walk,
i.e., for distinct \(z,w\), \(q(z,w)>0\) if and only if
\(w-z\in\mathrm{supp}(q_0)\setminus\{0\}\),
and \(Z\) is reversible with respect to a measure \(m\) satisfying
\(c_m\le m(z)\le C_m\) for all \(z\in\mathbb Z^d\), for some constants
\(0<c_m\le C_m<\infty\).
\end{enumerate}

Then there exists a constant \(C<\infty\) such that, for all
\(r\in\mathbb Z^d\) and all \(t\ge0\),
\begin{equation}
\label{eq:green-bound-all-d-abstract}
\sup_{z\in\mathbb Z^d}
\int_0^t
\mathbb P_z\!\big(Z(u)=r\big)\,du
\le
C
\begin{cases}
1+\sqrt t, & d=1,\\[0.3em]
1+\log(1+t), & d=2,\\[0.3em]
1, & d\ge3.
\end{cases}
\end{equation}
\end{proposition}

\begin{proof}
We denote by \(\mathcal E_m\) the Dirichlet form of \(Z\) in \(\ell^2(m)\),
which is given by
\[
\mathcal E_m(f,f)
:=
\frac12\sum_{z,w\in\mathbb Z^d}m(z)\,q(z,w)\bigl(f(w)-f(z)\bigr)^2 .
\]
We first compare it with the
form \(\mathcal E_0:=\mathcal E_0^{\mathbb Z^d}\) of the free walk. By
reversibility the conductances \(m(z)q(z,w)\) are symmetric in \((z,w)\).
For \(z\neq w\) with \(z\notin F\), we have
\(m(z)q(z,w)=m(z)q_0(w-z)\ge c_m\,q_0(w-z)\) by
assumption~\ref{hyp:green-finite-perturbation}, and by symmetry the same bound holds if \(w\notin F\). The
remaining pairs \(z,w\in F\) with \(w-z\in\mathrm{supp}(q_0)\setminus\{0\}\)
are finitely many, and for each of them \(m(z)q(z,w)>0\) by assumption~\ref{hyp:green-support-reversible}.
Hence there exists \(\kappa>0\) such that
\(m(z)q(z,w)\ge\kappa\,q_0(w-z)\) for all \(z\neq w\), i.e.,
\begin{equation}
\label{eq:energy-comparison}
\mathcal E_m(f,f)\ge\kappa\,\mathcal E_0(f,f).
\end{equation}

Combining \eqref{eq:energy-comparison} with
Lemma~\ref{lem:nash-free-walk} and with the bounds \(c_m\le m\le C_m\), we
obtain the following Nash inequality for the perturbed chain
\begin{equation}
\label{eq:nash-perturbed}
\begin{split}
\|f\|_{2,m}^{2+4/d}
&\le
C_m^{1+2/d}\,\|f\|_2^{2+4/d}
\\
&\le
C_m^{1+2/d}\,C_N\,\mathcal E_0(f,f)\,\|f\|_1^{4/d}
\\
&\le
\frac{C_m^{1+2/d}\,C_N}{\kappa\,c_m^{4/d}}\,
\mathcal E_m(f,f)\,\|f\|_{1,m}^{4/d} .
\end{split}
\end{equation}
Since the total jump rates are bounded, \(\mathcal E_m\) is a bounded
quadratic form on \(\ell^2(m)\), and \eqref{eq:nash-perturbed} extends by
truncation from finitely supported \(f\) to all
\(f\in\ell^1(m)\cap\ell^2(m)\).

Let \((P_t)_{t\ge0}\) be the semigroup of \(Z\), which is symmetric in
\(\ell^2(m)\). By \eqref{eq:nash-perturbed}
and~\cite[Theorem~(2.1)]{carlen1987upper}, with \(\nu=d\) and \(\delta=0\), there exists \(B<\infty\) such that
\(\|P_tf\|_\infty\le B\,t^{-d/2}\|f\|_{1,m}\) for all \(t>0\). We apply this
bound to \(f=\mathbf 1_{\{r\}}\), for which
\(P_tf(z)=\mathbb P_z(Z(t)=r)\) and \(\|f\|_{1,m}=m(r)\le C_m\). Since
probabilities are bounded by one, there exists a constant \(C>0\) such that
\begin{equation}
\label{eq:on-diagonal-abstract}
\mathbb P_z\bigl(Z(t)=r\bigr)\le \frac{C}{(1+t)^{d/2}},
\qquad z,r\in\mathbb Z^d,\ t\ge 0.
\end{equation}
Hence we have
\[
\int_0^t \mathbb P_z\bigl(Z(u)=r\bigr)\,du
\le
C\int_0^t (1+u)^{-d/2}\,du,
\qquad z\in\mathbb Z^d,\ t\ge0.
\]
Evaluating the right-hand side for \(d=1\), \(d=2\), and \(d\ge3\) (where the
integral converges) proves \eqref{eq:green-bound-all-d-abstract}.
\end{proof}

Hypothesis~\ref{hyp:green-support-reversible} does not cover the SEP(1) application, where the origin is
unreachable and the explicit reversible measure vanishes there. We therefore
also need the following estimate on a punctured lattice.

\begin{proposition}
\label{prop:green-bounds-punctured}
Let \(d\ge1\), \(v\in\mathbb Z^d\) and \(S_v:=\mathbb Z^d\setminus\{v\}\). Let
\(Z=(Z(t))_{t\ge0}\) be a continuous-time reversible Markov chain on \(S_v\)
with jump rates \(q(z,w)\), \(z,w\in S_v\), satisfying hypotheses~\ref{hyp:green-bounded-rates}--\ref{hyp:green-support-reversible} of
Proposition~\ref{prop:green-bounds-finite-perturb} with \(\mathbb Z^d\)
replaced by \(S_v\), where in~\ref{hyp:green-finite-perturbation} the finite set \(F\subset S_v\) contains
every \(z\in S_v\) with \(v-z\in\mathrm{supp}(q_0)\). Then the bound
\eqref{eq:green-bound-all-d-abstract} holds for all \(r\in S_v\), with the
supremum taken over \(z\in S_v\).
\end{proposition}

\begin{proof}
The proof of Proposition~\ref{prop:green-bounds-finite-perturb} applies with
\(\mathbb Z^d\) replaced by \(S_v\) throughout, using
Lemma~\ref{lem:nash-free-walk} with \(S=S_v\).
\end{proof}

\section{Proof of the Convergence Theorem}
\label{sec:proof-convergence}
Throughout the martingale-problem argument we use
\(C_c^\infty(\mathbb R^d)\) as the  class of test functions. For the tightness proof,
we also choose a countable family in this class that determines convergence,
together with the constant function \(1\).

\subsection{Dynkin martingales}
Let $\mathcal L^{(n),k}$ be the generator of $y^{(n),k}$.
For any bounded function $F$ on $\big(\frac{1}{n}\mathbb Z^d\big)^k$,
\[
M^F_t
:=F\big(y^{(n),k}(t)\big)-F\big(y^{(n),k}(0)\big)-\int_0^t (\mathcal L^{(n),k}F)\big(y^{(n),k}(s)\big)\,ds
\]
is a c\`adl\`ag martingale (with respect to the natural filtration).
Moreover, its predictable quadratic variation is given by
\begin{equation}
\label{eq:def-carre}
\langle M^F\rangle_t
=\int_0^t \Gamma^{(n),k}F\big(y^{(n),k}(s)\big)\,ds,
\qquad
\Gamma^{(n),k}F:=\mathcal L^{(n),k}(F^2)-2F\,\mathcal L^{(n),k}F.
\end{equation}
For \(\phi\in C_c^\infty(\mathbb R^d)\) and
\(y=(y_1,\dots,y_k)\in\big(\frac1n\mathbb Z^d\big)^k\) we consider the linear
observable
\begin{equation}
\label{eq:def-Fphi}
F_\phi(y):=\frac{1}{k}\sum_{i=1}^k \phi(y_i),
\end{equation}
so that \(F_\phi(y^{(n),k}(t))=\langle \mu_t^{(n),k},\phi\rangle\), and we denote
by \(M^{(n),\phi}:=M^{F_\phi}\) the associated Dynkin martingale, i.e.,
\begin{equation}
\label{eq:tightness-dynkin-Fphi}
M_t^{(n),\phi}
=
\langle \mu_t^{(n),k},\phi\rangle
-
\langle \mu_0^{(n),k},\phi\rangle
-
\int_0^t
(\mathcal L^{(n),k}F_\phi)(y^{(n),k}(s))\,ds .
\end{equation}

\subsection{Tightness}
\label{subsec:tightness}

Tightness of the empirical measure process \(\mu^{(n),k}\) in
\(D([0,T],\mathcal M_1(\mathbb R^d))\) follows from two estimates: compact
containment, which we obtain from the second moment of \(\mu_t^{(n),k}\), and
tightness of the real-valued projections \(\langle\mu_t^{(n),k},\phi\rangle\),
which we obtain from bounds on the drift and on the carr\'e du champ. The projection criterion of Roelly-Coppoletta
\cite[Thm.~2.1]{roelly1986criterion} then combines the two.

\begin{proposition}
\label{prop:compact-containment}
Assume that
\[
y^{(n),k}(0)\to y^k=(y_1,\dots,y_k)
\qquad\text{in }(\mathbb R^d)^k.
\]
Then, for every \(T>0\),
\[
\lim_{A\to\infty}
\sup_{n\ge1}
\mathbb P\left(
\sup_{0\le t\le T}
\int_{\mathbb R^d}|x|^2\,\mu_t^{(n),k}(dx)>A
\right)
=0.
\]
Consequently, for every \(\varepsilon>0\) and \(T>0\), there exists a compact set
\(\mathcal K_\varepsilon\subset \mathcal M_1(\mathbb R^d)\) such that
\[
\inf_{n\ge1}
\mathbb P\left(
\mu_t^{(n),k}\in \mathcal K_\varepsilon
\text{ for all }t\in[0,T]
\right)
\ge 1-\varepsilon .
\]
\end{proposition}

\begin{proof}
Define
\[
V(y)
:=
\frac1k\sum_{i=1}^k |y_i|^2,
\qquad
y=(y_1,\dots,y_k)\in\Big(\frac1n\mathbb Z^d\Big)^k,
\]
so that \(V(y^{(n),k}(t))=\int_{\mathbb R^d}|x|^2\,\mu_t^{(n),k}(dx)\). A jump of
particle \(i\) by \(r/n\) changes \(V\) by
\(\frac1k\bigl(2y_i\cdot \frac rn+\frac{|r|^2}{n^2}\bigr)\), which on the event
\(y_j=y_i+r/n\) equals \(\frac1k(|y_j|^2-|y_i|^2)\). In the free part of the
generator the first-order term vanishes because \(p\) is symmetric, and in the
interaction part the \((i,j)\) summand cancels exactly against the \((j,i)\)
summand, again by the symmetry of \(p\). Hence, reading
\(p\bigl(n(y_j-y_i)\bigr)\) as zero unless \(n(y_j-y_i)\in\mathbb Z^d\), we
obtain
\begin{equation}
\label{eq:LV-second-moment}
(\mathcal L^{(n),k} \, V)(y)
=
\alpha\,\mathrm{tr}(a)
+
\frac{\sigma n^2}{k}
\sum_{i\neq j}
p\bigl(n(y_j-y_i)\bigr)
\bigl(|y_j|^2-|y_i|^2\bigr)
=
\alpha\,\mathrm{tr}(a).
\end{equation}

For \(A>0\) let \(\tau_A:=\inf\{t\ge0:\,V(y^{(n),k}(t))\ge A\}\). For fixed \(n\)
and \(A\) the process stopped at \(\tau_A\) takes values in a bounded set,
since \(p\) has finite range, hence Dynkin's formula applies to \(V\) up to
time \(\tau_A\). Since \(V\) is nonnegative, by \eqref{eq:LV-second-moment} we
obtain
\begin{align*}
A\,\mathbb P(\tau_A\le T)
&\le
\mathbb E\bigl[V\bigl(y^{(n),k}(T\wedge\tau_A)\bigr)\bigr]
\\
&=
V\bigl(y^{(n),k}(0)\bigr)+\alpha\,\mathrm{tr}(a)\,\mathbb E[T\wedge\tau_A]
\\
&\le
C_0+\alpha\,\mathrm{tr}(a)\,T ,
\end{align*}
where \(C_0:=\sup_nV(y^{(n),k}(0))<\infty\) because \(y^{(n),k}(0)\to y^k\).
If \(V(y^{(n),k}(t))>A\) for some \(t\le T\), then \(\tau_A\le T\), hence
letting \(A\to\infty\) proves the first claim.

For the compact-containment statement, the sublevel set
\(\mathcal K_A:=\{\nu\in \mathcal M_1(\mathbb R^d):\int_{\mathbb R^d}|x|^2\,\nu(dx)\le A\}\)
is compact in the weak topology: it is tight by Markov's inequality
(\(\nu(\{|x|>R\})\le A/R^2\) uniformly on \(\mathcal K_A\)), hence relatively
compact by Prokhorov's theorem, and closed because \(x\mapsto|x|^2\) is
nonnegative and lower semicontinuous. Choosing \(A\) such that
\(C_0+\alpha\,\mathrm{tr}(a)T\le\varepsilon A\)
gives the compact-containment condition.
\end{proof}

We next prove tightness of the projections against a convergence-determining
class of smooth compactly supported test functions.

\begin{proposition}
\label{prop:tightness-projections}
For every \(\phi\in C_c^\infty(\mathbb R^d)\), the sequence of real-valued processes
\[
\left(
\langle \mu_t^{(n),k},\phi\rangle
\right)_{0\le t\le T}
\]
is tight in \(D([0,T],\mathbb R)\). Moreover, every limit point is supported on
\(C([0,T],\mathbb R)\).
\end{proposition}

\begin{proof}
Fix \(\phi\in C_c^\infty(\mathbb R^d)\). By \eqref{eq:def-carre}, both
\(M^{(n),\phi}\) and
\(\bigl(M_t^{(n),\phi}\bigr)^2-\int_0^t\Gamma^{(n),k}F_\phi(y^{(n),k}(s))\,ds\)
are martingales.
We verify the Aldous-Rebolledo tightness
criterion~\cite{aldous1978stopping,rebolledo1980existence} for this real-valued
semimartingale, in the form used in
Roelly-Coppoletta~\cite[Thm.~2.3]{roelly1986criterion}: it suffices to control
the second moments and the increments of the finite-variation part and of the
martingale bracket along bounded stopping times. The second moments are
immediate, \(\sup_{n}\sup_{t\le T}\mathbb E[|\langle
\mu_t^{(n),k},\phi\rangle|^2]\le\|\phi\|_\infty^2\).

For the drift, Proposition~\ref{prop:drift_closure} gives
\[
(\mathcal L^{(n),k}F_\phi)(y)
=\langle \mu_y,\mathcal L^{(n),1}\phi\rangle,
\qquad
(\mathcal L^{(n),1}\phi)(x)
=n^2\alpha\sum_r p(r)\bigl[\phi(x+\tfrac rn)-\phi(x)\bigr].
\]
In the second-order Taylor expansion of \(\phi(x+\frac rn)-\phi(x)\) the
first-order term vanishes after summation against the symmetric \(p\). Hence,
writing \(\|D^2\phi\|_\infty:=\sup_x\|D^2\phi(x)\|_{\mathrm{op}}\) for the
supremum of the operator norm of the Hessian, we obtain
\begin{equation}
\label{eq:gamma1-bound-tightness}
\sup_{n\ge1}\sup_y
\left|
(\mathcal L^{(n),k}F_\phi)(y)
\right|
\le
\frac{\alpha}{2}\,\mathrm{tr}(a)\,\|D^2\phi\|_\infty .
\end{equation}

For the carr\'e du champ, a jump of particle \(i\) by \(r/n\) changes
\(F_\phi\) by at most \(\frac{|r|}{kn}\|\nabla\phi\|_\infty\). Moreover,
uniformly in \(y\), the rates are nonnegative and bounded by \(\alpha+k-1\)
(for \(\sigma=-1\) they are bounded by \(\alpha\)). Hence we obtain
\begin{equation}
\label{eq:gamma2-bound-tightness}
\sup_{n\ge1}\sup_y
\Gamma^{(n),k}F_\phi(y)
\le
\frac{\alpha+k-1}{k}\,\mathrm{tr}(a)\,\|\nabla\phi\|_\infty^2 .
\end{equation}
In what follows \(C_\phi\) denotes a constant larger than the right-hand sides
of \eqref{eq:gamma1-bound-tightness} and \eqref{eq:gamma2-bound-tightness}.

Let \(\tau_n\le T\) be stopping times and let \(\theta_n\downarrow0\). By
\eqref{eq:gamma1-bound-tightness}, the drift part of the increment of
\(\langle\mu^{(n),k},\phi\rangle\) over \([\tau_n,\tau_n+\theta_n]\) is bounded
by \(C_\phi\theta_n\), and by optional stopping together with
\eqref{eq:gamma2-bound-tightness} we have
\[
\mathbb E\Bigl[\bigl|M_{\tau_n+\theta_n}^{(n),\phi}-M_{\tau_n}^{(n),\phi}\bigr|^2\Bigr]
\le C_\phi\theta_n .
\]
Hence
\(\langle \mu_{\tau_n+\theta_n}^{(n),k},\phi\rangle-\langle \mu_{\tau_n}^{(n),k},\phi\rangle\to0\)
in probability, which is Aldous' tightness condition for the real-valued
projections.

Finally, jumps vanish: at a microscopic jump exactly one particle moves by
\(r/n\) with \(r\in\mathrm{supp}(p)\), so
\(|F_\phi(y^{i,i+r/n})-F_\phi(y)|\le\|\nabla\phi\|_\infty|r|/(kn)\) and
\[
\sup_{0\le t\le T}
\left|
\langle \mu_t^{(n),k},\phi\rangle
-
\langle \mu_{t-}^{(n),k},\phi\rangle
\right|
\le
\frac{C_\phi}{n}
\longrightarrow 0.
\]
Thus the projections are tight in \(D([0,T],\mathbb R)\), and every limit point
is supported on \(C([0,T],\mathbb R)\).
\end{proof}

\begin{theorem}
\label{thm:tightness-empirical-measures}
Assume that
\[
y^{(n),k}(0)\to y^k
\qquad\text{in }(\mathbb R^d)^k.
\]
Then, for every \(T>0\), the sequence
\[
\bigl(\mu^{(n),k}\bigr)_{n\ge1}
\]
is tight in
\[
D([0,T],\mathcal M_1(\mathbb R^d)).
\]
Moreover, every limit point is supported on \(C([0,T],\mathcal M_1(\mathbb R^d))\).
\end{theorem}

\begin{proof}
We apply the projection criterion of
Roelly-Coppoletta~\cite[Thm.~2.1]{roelly1986criterion}. Choose a countable
family \(\mathcal D:=\{\phi_m:m\ge1\}\subset C_c^\infty(\mathbb R^d)\), dense in
\(C_0(\mathbb R^d)\) for the uniform norm, and augment it by the constant
function \(\phi_0:=1\). For \(m\ge1\),
Proposition~\ref{prop:tightness-projections} gives tightness of
\((\langle \mu_t^{(n),k},\phi_m\rangle)_{0\le t\le T}\) in \(D([0,T],\mathbb R)\),
while for \(m=0\) the projection is identically one and hence trivially tight.
All projected laws required by the criterion are therefore tight, and
\((\mu^{(n),k})_{n\ge1}\) is tight in
\(D([0,T],\mathcal M_F(\mathbb R^d))\). Since each \(\mu_t^{(n),k}\) is a
probability measure and the total-mass projection is identically one, every
subsequential limit is \(\mathcal M_1(\mathbb R^d)\)-valued; together with the
compact-containment estimate of Proposition~\ref{prop:compact-containment} this
gives tightness in \(D([0,T],\mathcal M_1(\mathbb R^d))\) for the weak topology.

For the continuity of limit points, we metrise weak convergence
on \(\mathcal M_1(\mathbb R^d)\) by
\(d_{\mathcal M_1}(\nu,\rho):=\sum_{m\ge1}2^{-m}(|\langle\nu-\rho,\phi_m\rangle|\wedge1)\),
with \((\phi_m)_{m\ge1}\) bounded and convergence determining. By the last
estimate in the proof of Proposition~\ref{prop:tightness-projections}, for
every \(M\ge1\) we have
\[
\sup_{0\le t\le T}
d_{\mathcal M_1}\bigl(\mu_t^{(n),k},\mu_{t-}^{(n),k}\bigr)
\le
\sum_{m=1}^M2^{-m}\Bigl(\frac{C_{\phi_m}}{n}\wedge1\Bigr)+2^{-M} .
\]
Letting \(n\to\infty\) and then \(M\to\infty\), the maximal jump size vanishes,
and together with tightness this implies that every
subsequential limit is supported on continuous
\(\mathcal M_1(\mathbb R^d)\)-valued paths.
\end{proof}

\subsection{Closing the drift}

Let $\Psi$ be the projection from coordinates to configurations defined
in~\eqref{eq:def-Psi-main}, and denote by $\eta_y:=\Psi(y)$ the configuration
associated with $y$. Then the linear observable $F_\phi$
of~\eqref{eq:def-Fphi} admits the configuration representation
\begin{equation}
\label{eq:Fphi-conf}
F_\phi(y)
=
\frac1k \sum_{x\in\frac1n\mathbb Z^d}\eta_y(x)\phi(x)
=
\frac1k \sum_{x\in\frac1n\mathbb Z^d} D(\delta_x,\eta_y)\phi(x),
\end{equation}
where $D$ denotes the one-particle self-duality function.

\begin{proposition}
\label{prop:drift_closure}
For every $\phi\in C_c^\infty(\mathbb R^d)$,
\begin{equation}
\label{eq:closed-drift}
(\mathcal L^{(n),k}F_\phi)(y)
=
F_{\mathcal L^{(n),1}\phi}(y)
\qquad
\text{for all }y\in\Big(\frac1n\mathbb Z^d\Big)^k .
\end{equation}
\end{proposition}

\begin{proof}
Using \eqref{eq:Fphi-conf} we obtain
\begin{equation}
\label{eq:drift_start}
(\mathcal L^{(n),k}F_\phi)(y)
=
\frac1k
\sum_{x\in\frac1n\mathbb Z^d}
(\mathcal L^{(n),k} D(\delta_x,\eta_y))\,
\phi(x).
\end{equation}

Since the coordinate and configuration processes are compatible (Lemma~2.5 in~\cite{carinci2021consistent}),
the coordinate generator
$\mathcal L^{(n),k}$ acting on $D(\delta_x,\eta_y)$ coincides with the
configuration generator acting on the configuration variable $\eta_y$.
Therefore
\[
(\mathcal L^{(n),k} D(\delta_x,\eta_y))
=
(\mathcal L^{(n)} D(\delta_x,\cdot))(\eta_y).
\]

By the one-particle self-duality relation
\eqref{eq:occupation-duality-generator},
\[
(\mathcal L^{(n)} D(\delta_x,\cdot))(\eta)
=
(\mathcal L^{(n),1} D(\delta_{\cdot},\eta))(x),
\]
so that \eqref{eq:drift_start} becomes
\[
(\mathcal L^{(n),k}F_\phi)(y)
=
\frac1k
\sum_{x\in\frac1n\mathbb Z^d}
(\mathcal L^{(n),1} D(\delta_{\cdot},\eta_y))(x)
\,\phi(x).
\]

Finally, since $\mathcal L^{(n),1}$ is reversible with respect to the
counting measure on $\frac1n\mathbb Z^d$, self-adjointness yields
\[
\sum_x (\mathcal L^{(n),1} D(\delta_{\cdot},\eta_y))(x)\,\phi(x)
=
\sum_x D(\delta_x,\eta_y)\,(\mathcal L^{(n),1}\phi)(x).
\]

Substituting this into the previous expression and using
\eqref{eq:Fphi-conf} gives
\[
(\mathcal L^{(n),k}F_\phi)(y)
=
\frac1k
\sum_x D(\delta_x,\eta_y)\,
(\mathcal L^{(n),1}\phi)(x)
=
F_{\mathcal L^{(n),1}\phi}(y),
\]
which proves \eqref{eq:closed-drift}.
\end{proof}

\subsection{Closing the quadratic variation}

Fix $\phi\in C_c^\infty(\mathbb R^d)$ and let $F_\phi$ be the linear observable
defined in \eqref{eq:def-Fphi}. For $r\in\mathbb Z^d$ and
$x\in \frac1n\mathbb Z^d$, write
\[
\Delta_r^n \phi(x):=\phi\!\left(x+\frac{r}{n}\right)-\phi(x),
\qquad
q_n\phi(x):=
n^2\sum_{r\in\mathbb Z^d} p(r)\,\big(\Delta_r^n\phi(x)\big)^2.
\]

Recall that
\[
\mu_y:=\frac1k\sum_{i=1}^k\delta_{y_i},
\qquad y=(y_1,\dots,y_k)\in \Big(\frac1n\mathbb Z^d\Big)^k.
\]

\begin{proposition}
\label{prop:qv-decomposition}
For every $y\in \big(\frac1n\mathbb Z^d\big)^k$,
\begin{equation}
\label{eq:qv-decomposition}
\Gamma^{(n),k}F_\phi(y)
=
\frac{\alpha}{k}\, F_{q_n\phi}(y)
+
R_{\mathrm{int}}^{(n),k}(y;\phi),
\end{equation}
where
\begin{equation}
\label{eq:Rint-def}
R_{\mathrm{int}}^{(n),k}(y;\phi)
:=
\frac{\sigma\,n^2}{k^2}
\sum_{r\in\mathbb Z^d} p(r)
\sum_{\substack{i,j=1\\ i\neq j}}^k
\mathbf 1_{\{y_j=y_i+\frac{r}{n}\}}
\big(\Delta_r^n\phi(y_i)\big)^2.
\end{equation}
Moreover, there exists a constant $C_\phi<\infty$ such that
\begin{equation}
\label{eq:qn-to-grad}
\sup_{x\in \frac1n\mathbb Z^d}
\left|
q_n\phi(x)-\nabla\phi(x)\cdot a\nabla\phi(x)
\right|
\le \frac{C_\phi}{n}.
\end{equation}
\end{proposition}

\begin{proof}
By definition of the generator $\mathcal L^{(n),k}$, a jump of particle $i$
by $r/n$ changes $F_\phi$ by
\[
F_\phi(y^{i,i+\frac{r}{n}})-F_\phi(y)
=
\frac1k\,\Delta_r^n\phi(y_i).
\]
Hence, using the definition \eqref{eq:def-carre} of the carr\'e du champ,
\begin{align*}
\Gamma^{(n),k}F_\phi(y)
&=
n^2\sum_{i=1}^k\sum_{r\in\mathbb Z^d}
p(r)\Big(\alpha+\sigma\sum_{j\neq i}\mathbf 1_{\{y_j=y_i+\frac{r}{n}\}}\Big)
\left(\frac1k\Delta_r^n\phi(y_i)\right)^2 \\
&=
\frac{\alpha n^2}{k^2}
\sum_{i=1}^k\sum_{r\in\mathbb Z^d}
p(r)\big(\Delta_r^n\phi(y_i)\big)^2
+
R_{\mathrm{int}}^{(n),k}(y;\phi).
\end{align*}
Since
\[
\frac1k\sum_{i=1}^k q_n\phi(y_i)= F_{q_n\phi}(y),
\]
this proves \eqref{eq:qv-decomposition}.

\medskip
For \eqref{eq:qn-to-grad}, since $\phi\in C_c^\infty(\mathbb R^d)$, Taylor's formula yields
\[
\Delta_r^n\phi(x)
=
\frac1n\,r\cdot \nabla\phi(x)+R_n(x,r),
\qquad
|R_n(x,r)|
\le
\frac12 \|D^2\phi\|_\infty \frac{|r|^2}{n^2},
\]
uniformly in $x$. Hence
\begin{align*}
\left|
n^2\big(\Delta_r^n\phi(x)\big)^2-(r\cdot \nabla\phi(x))^2
\right|
&\le
2n\,|r\cdot \nabla\phi(x)|\,|R_n(x,r)|
+
n^2 |R_n(x,r)|^2 \\
&\le
\|\nabla\phi\|_\infty \|D^2\phi\|_\infty \frac{|r|^3}{n}
+
\frac14 \|D^2\phi\|_\infty^2 \frac{|r|^4}{n^2}.
\end{align*}
Since $p$ has finite range, summation over $r$ gives
\[
\sup_{x\in\frac1n\mathbb Z^d}
\left|
q_n\phi(x)-\sum_{r\in\mathbb Z^d}p(r)(r\cdot \nabla\phi(x))^2
\right|
\le \frac{C(\phi,p)}{n}.
\]
Finally,
\[
\sum_{r\in\mathbb Z^d}p(r)(r\cdot \nabla\phi(x))^2
=
\nabla\phi(x)\cdot a\nabla\phi(x),
\]
which proves \eqref{eq:qn-to-grad}.
\end{proof}

\begin{remark}
The second term
$R_{\mathrm{int}}^{(n),k}$ is supported only on microscopic near-collisions
$y_j=y_i+\frac{r}{n}$ and therefore carries the interaction correction.
\end{remark}

The following proposition uses Corollary~\ref{cor:time_integrated_pair_reduction}
to reduce the time-integrated interaction term~\eqref{eq:Rint-def} to estimates
on two-particle systems started from the corresponding initial pairs.

\begin{proposition}
\label{prop:qv-two-particle-reduction}
Let
\[
\mathscr R_T^{(n),k}(\phi)
:=
\int_0^T
\big|R_{\mathrm{int}}^{(n),k}(y^{(n),k}(s);\phi)\big|\,ds .
\]
Write
\[
x^{(n),k}(0)=(x^{(n)}_1,\dots,x^{(n)}_k)\in(\mathbb Z^d)^k
\]
for the unscaled initial configuration corresponding to
$y^{(n),k}(0)=x^{(n),k}(0)/n$. Then, for every $T>0$,
\begin{align}
\mathbb E\big[\mathscr R_T^{(n),k}(\phi)\big]
&\le
\frac{|\sigma|}{k^2}
\sum_{r\in\mathbb Z^d}p(r)\,C_{\phi,r}
\sum_{\substack{1\le i,j\le k\\ i\neq j}}
\int_0^T
\mathbb P^{(2)}_{(x_i^{(n)},x_j^{(n)})}
\left(
y^{(n),2}_2(s)=y^{(n),2}_1(s)+\frac{r}{n}
\right)\,ds
\label{eq:two-particle-reduction-pairwise}
\\
&=
\frac{|\sigma|}{k^2 n^2}
\sum_{r\in\mathbb Z^d}p(r)\,C_{\phi,r}
\sum_{\substack{1\le i,j\le k\\ i\neq j}}
\int_0^{n^2T}
\mathbb P^{(2)}_{(x_i^{(n)},x_j^{(n)})}
\left(
x^{(2)}_2(u)=x^{(2)}_1(u)+r
\right)\,du .
\label{eq:two-particle-reduction-unscaled-pairwise}
\end{align}
Equivalently, the rougher uniform bound
\begin{align}
\mathbb E\big[\mathscr R_T^{(n),k}(\phi)\big]
&\le
\frac{|\sigma|(k-1)}{k n^2}
\sum_{r\in\mathbb Z^d}p(r)\,C_{\phi,r}
\sup_{z\in\mathbb Z^d}
\int_0^{n^2T}
\mathbb P_z^{\Xi}
\big(\Xi(u)=r\big)\,du
\label{eq:two-particle-reduction-uniform}
\end{align}
holds, where $\Xi(u):=x^{(2)}_2(u)-x^{(2)}_1(u)$ denotes the two-particle
difference process started from $z=x_j^{(n)}-x_i^{(n)}$.

Here
\[
C_{\phi,r}:=
\sup_{n\ge1}\sup_{x\in\mathbb R^d}
n^2\big(\Delta_r^n\phi(x)\big)^2
\le |r|^2\|\nabla\phi\|_\infty^2 <\infty,
\]
by the mean-value theorem.
\end{proposition}

\begin{proof}
Bounding $n^2\big(\Delta_r^n\phi(y_i)\big)^2$ by $C_{\phi,r}$ in \eqref{eq:Rint-def} gives
\[
\big|R_{\mathrm{int}}^{(n),k}(y;\phi)\big|
\le
\frac{|\sigma|}{k^2}
\sum_{r\in\mathbb Z^d}p(r)\,C_{\phi,r}
\sum_{\substack{1\le i,j\le k\\ i\neq j}}
\mathbf 1_{\{y_j=y_i+r/n\}} .
\]
Integrating in time and taking expectations gives
\[
\mathbb E\big[\mathscr R_T^{(n),k}(\phi)\big]
\le
\frac{|\sigma|}{k^2}
\sum_{r\in\mathbb Z^d}p(r)\,C_{\phi,r}
\sum_{\substack{1\le i,j\le k\\ i\neq j}}
\int_0^T
\mathbb P^{(k)}_{x^{(n),k}(0)}
\left(
y^{(n),k}_j(s)=y^{(n),k}_i(s)+\frac r n
\right)\,ds .
\]
Applying Corollary~\ref{cor:time_integrated_pair_reduction} with
\[
h_r(u,v):=\mathbf 1_{\{v=u+r/n\}}
\]
yields the pairwise representation
\eqref{eq:two-particle-reduction-pairwise}. The change of variables
$u=n^2s$ gives \eqref{eq:two-particle-reduction-unscaled-pairwise}.
Finally, since the event
\[
x^{(2)}_2(u)=x^{(2)}_1(u)+r
\]
is the same as $\Xi(u)=r$, and since there are $k(k-1)$ ordered pairs, taking
the supremum over the initial difference
$z=x_j^{(n)}-x_i^{(n)}$ gives
\eqref{eq:two-particle-reduction-uniform}.
\end{proof}

\begin{proposition}
\label{prop:vanish-exp-loc-time}
Let
\begin{equation}
\label{eq:def-diff-proc}
\Xi(u):=x^{(2)}_2(u)-x^{(2)}_1(u)
\end{equation}
denote the two-particle difference process. Then, for every \(T>0\) and every
\(r\) in the finite support of \(p\),
\begin{equation}
\label{eq:two-particle-occupation-vanishing}
\frac1{n^2}
\int_0^{n^2T}
\mathbb P\!\left(
x^{(2)}_2(u)=x^{(2)}_1(u)+r
\right)\,du
=
\frac1{n^2}
\int_0^{n^2T}
\mathbb P\!\left(\Xi(u)=r\right)\,du
\longrightarrow 0,
\qquad n\to\infty.
\end{equation}
More precisely,
\begin{equation}
\label{eq:occupation-rates-all-d}
\frac1{n^2}
\int_0^{n^2T}
\mathbb P\!\left(\Xi(u)=r\right)\,du
=
\begin{cases}
O(n^{-1}), & d=1,\\[0.4em]
O\!\left(\dfrac{\log n}{n^2}\right), & d=2,\\[0.9em]
O(n^{-2}), & d\ge 3.
\end{cases}
\end{equation}
\end{proposition}

\begin{proof}
By definition of \(\Xi\),
\[
\mathbb P\!\left(x^{(2)}_2(u)=x^{(2)}_1(u)+r\right)
=
\mathbb P\!\left(\Xi(u)=r\right),
\]
so it suffices to estimate the truncated Green function of \(\Xi\).

We first compute the generator of \(\Xi\). Let \(f\colon \mathbb Z^d\to\mathbb R\) be
bounded and define
\[
F(x_1,x_2):=f(x_2-x_1).
\]
Writing \(z=x_2-x_1\), the two-particle generator obtained from
\eqref{eq:unified-k-generator} gives
\[
(L^{(2)}F)(x_1,x_2)
=
\sum_{\rho\in\mathbb Z^d} p(\rho)\,
\bigl(\alpha+\sigma \mathbf 1_{\{x_2=x_1+\rho\}}\bigr)\,
\bigl(f(z-\rho)-f(z)\bigr)
\]
\[
\hspace{5em}
+
\sum_{\rho\in\mathbb Z^d} p(\rho)\,
\bigl(\alpha+\sigma \mathbf 1_{\{x_1=x_2+\rho\}}\bigr)\,
\bigl(f(z+\rho)-f(z)\bigr).
\]
Since
\[
\mathbf 1_{\{x_2=x_1+\rho\}}=\mathbf 1_{\{z=\rho\}},
\qquad
\mathbf 1_{\{x_1=x_2+\rho\}}=\mathbf 1_{\{z=-\rho\}},
\]
we obtain
\begin{equation}
\label{eq:generator-relative-coordinate}
\begin{aligned}
(L_\Xi f)(z)
&=
\sum_{\rho\in\mathbb Z^d} p(\rho)\,
\bigl(\alpha+\sigma \mathbf 1_{\{z=\rho\}}\bigr)\,
\bigl(f(z-\rho)-f(z)\bigr)
\\
&\quad+
\sum_{\rho\in\mathbb Z^d} p(\rho)\,
\bigl(\alpha+\sigma \mathbf 1_{\{z=-\rho\}}\bigr)\,
\bigl(f(z+\rho)-f(z)\bigr).
\end{aligned}
\end{equation}

This representation shows that \(\Xi\) has bounded total jump rates and finite
jump range. Moreover, the perturbation terms involving the indicators can be
nonzero only when
\[
z\in F:=\operatorname{supp}(p)\cup(-\operatorname{supp}(p)).
\]
Since \(p\) is symmetric, \(F=\operatorname{supp}(p)\), which is finite.

For \(z\notin F\), the indicator terms vanish and
\eqref{eq:generator-relative-coordinate} reduces to
\[
(L_\Xi f)(z)
=
\alpha \sum_{\rho\in\mathbb Z^d} p(\rho)\bigl(f(z-\rho)-f(z)\bigr)
+
\alpha \sum_{\rho\in\mathbb Z^d} p(\rho)\bigl(f(z+\rho)-f(z)\bigr).
\]
Using the symmetry of \(p\), this can be rewritten as
\[
(L_\Xi f)(z)
=
\sum_{b\in\mathbb Z^d} q_0(b)\bigl(f(z+b)-f(z)\bigr),
\qquad z\notin F,
\]
where
\[
q_0(b):=2\alpha\,p(b).
\]
Hence, outside the finite defect set \(F\), the process \(\Xi\) coincides with a
symmetric finite-range random walk on \(\mathbb Z^d\). Since
\(\operatorname{supp}(q_0)=\operatorname{supp}(p)\) and \(\operatorname{supp}(p)\)
generates \(\mathbb Z^d\), this underlying walk is irreducible.

The process \(\Xi\) is reversible. Indeed, if
\[
m(z):=\alpha+\sigma\,\mathbf 1_{\{z=0\}},
\]
then, by the symmetry of \(p\), the jump rates in
\eqref{eq:generator-relative-coordinate} are given by
\[
q(z,w)=2\,p(w-z)\,m(w),
\qquad z\neq w,
\]
and hence they satisfy the detailed balance relation
\[
m(z)\,q(z,w)=2\,p(w-z)\,m(z)\,m(w)=m(w)\,q(w,z),
\qquad z\neq w.
\]
If \((\sigma,\alpha)\neq(-1,1)\), then \(m(0)=\alpha+\sigma>0\), \(m\) is
bounded above and below by positive constants, and all the hypotheses of
Proposition~\ref{prop:green-bounds-finite-perturb} are satisfied by
\(Z=\Xi\) outright. If \((\sigma,\alpha)=(-1,1)\) (the classical simple
exclusion process, capacity one per site), then \(m(0)=0\): the origin is
unreachable for \(\Xi\), consistently with the fact that the labelled SEP(1)
process never lets two coordinates collide, so that \(\Xi\) started away
from \(0\) in fact never visits it (cf.\
\(E_k^{\mathrm{SEP}}\) in Section~\ref{sec:microscopic-model}). In this case
\(\Xi\) is a Markov chain on \(\mathbb Z^d\setminus\{0\}\), where
\(m\equiv1\), and it satisfies the hypotheses of
Proposition~\ref{prop:green-bounds-punctured} with \(v=0\) and with the finite
set \(F\setminus\{0\}\). Hence the bound
\eqref{eq:green-bound-all-d-abstract} holds for \(\Xi\) for
every \(r\neq0\). We only ever invoke it below at \(r\in\operatorname{supp}(p)\);
if \(0\in\operatorname{supp}(p)\) the case \(r=0\) is trivial in any case,
since a displacement \(r=0\) leaves the pair unchanged and contributes
nothing to the interaction term \(R_{\mathrm{int}}^{(n),k}\) in
\eqref{eq:Rint-def}. In all cases, for every \(r\in\operatorname{supp}(p)\),
the bound \eqref{eq:green-bound-all-d-abstract} with \(t=n^2T\), divided by
\(n^2\), gives
\[
\frac1{n^2}
\int_0^{n^2T}
\mathbb P\!\left(\Xi(u)=r\right)\,du
\le
\frac{C}{n^2}
\begin{cases}
1+n\sqrt T, & d=1,\\[0.3em]
1+\log(1+n^2T), & d=2,\\[0.3em]
1, & d\ge3,
\end{cases}
\]
which proves \eqref{eq:occupation-rates-all-d}. In particular,
\eqref{eq:two-particle-occupation-vanishing} follows.
\end{proof}

As a consequence we have that the interaction component vanishes in the $L^1$ sense.

\begin{corollary}
\label{cor:qv-interaction-vanishes}
For every $\phi\in C_c^\infty(\mathbb R^d)$,
\begin{equation}
\label{eq:Rint-vanishes}
\mathbb E\!\left[
\int_0^T \big|R_{\mathrm{int}}^{(n),k}(y^{(n),k}(s);\phi)\big|\,ds
\right]
\longrightarrow 0,
\qquad n\to\infty.
\end{equation}
Consequently,
\begin{equation}
\label{eq:qv-closes-final}
\mathbb E\!\left[
\sup_{t\in[0,T]}
\left|
\langle M^{F_\phi}\rangle_t
-
\frac{\alpha}{k}\int_0^t
\langle \mu_s^{(n),k},\nabla\phi\cdot a\nabla\phi\rangle\,ds
\right|
\right]
\longrightarrow 0.
\end{equation}
\end{corollary}

\begin{proof}
By Proposition~\ref{prop:qv-two-particle-reduction} and the finite-range assumption
on $p$, \eqref{eq:two-particle-occupation-vanishing} implies \eqref{eq:Rint-vanishes}.

Moreover, by Proposition~\ref{prop:qv-decomposition},
\[
\langle M^{F_\phi}\rangle_t
=
\int_0^t \Gamma^{(n),k}F_\phi(y^{(n),k}(s))\,ds
=
\frac{\alpha}{k}\int_0^t \langle \mu_s^{(n),k},q_n\phi\rangle\,ds
+
\int_0^t R_{\mathrm{int}}^{(n),k}(y^{(n),k}(s);\phi)\,ds.
\]
Therefore
\begin{align*}
&\sup_{t\in[0,T]}
\left|
\langle M^{F_\phi}\rangle_t
-
\frac{\alpha}{k}\int_0^t
\langle \mu_s^{(n),k},\nabla\phi\cdot a\nabla\phi\rangle\,ds
\right| \\
&\qquad\le
\frac{\alpha}{k}\int_0^T
\Big|
\langle \mu_s^{(n),k},q_n\phi-\nabla\phi\cdot a\nabla\phi\rangle
\Big|\,ds
+
\int_0^T
\big|R_{\mathrm{int}}^{(n),k}(y^{(n),k}(s);\phi)\big|\,ds.
\end{align*}
Since $\mu_s^{(n),k}$ is a probability measure and
\eqref{eq:qn-to-grad} holds uniformly,
\[
\Big|
\langle \mu_s^{(n),k},q_n\phi-\nabla\phi\cdot a\nabla\phi\rangle
\Big|
\le \frac{C_\phi}{n}.
\]
Taking expectations and using \eqref{eq:Rint-vanishes} proves
\eqref{eq:qv-closes-final}.
\end{proof}

\subsection{Identification of limit points}
\label{subsec:identification-limit-points}

We now identify the subsequential limits provided by tightness. The goal is to show
that every limit point satisfies the Dean-Kawasaki martingale problem. Uniqueness of
that martingale problem then identifies the limit and removes the dependence on the
chosen subsequence.

\begin{proposition}
\label{prop:identification-limit-points}
Assume the hypotheses of Theorem~\ref{thm:empirical-to-dk}. Let
\((\mu^{(n_j),k})_{j\ge1}\) be a subsequence such that
\[
\mu^{(n_j),k}\Longrightarrow \mu
\qquad\text{in }D([0,T],\mathcal M_1(\mathbb R^d)).
\]
Then \(\mu\) solves the Dean-Kawasaki martingale problem of
Definition~\ref{def:DK-mp} with parameter \(k\) and initial condition
\(\mu^k_0\).
\end{proposition}

\begin{proof}
Fix \(\phi\in C_c^\infty(\mathbb R^d)\), and relabel the subsequence \(n_j\) as
\(n\). Let \(\mathcal A\phi:=\frac{\alpha}{2}\,a:D^2\phi\) and
\[
\varepsilon_n
:=
\sup_{x\in \frac1n\mathbb Z^d}
\left|
\mathcal L^{(n),1}\phi(x)-\mathcal A\phi(x)
\right| ,
\]
which tends to zero by Taylor expansion, since \(p\) is symmetric and has
finite range. By Proposition~\ref{prop:drift_closure}, the Dynkin martingale
\(M^{(n),\phi}\) of \eqref{eq:tightness-dynkin-Fphi} has drift
\(\langle\mu_s^{(n),k},\mathcal L^{(n),1}\phi\rangle\). Hence the process
\[
\widetilde M^{(n),\phi}_t
:=
\langle \mu^{(n),k}_t,\phi\rangle
-
\langle \mu^{(n),k}_0,\phi\rangle
-
\int_0^t
\langle \mu^{(n),k}_s,\mathcal A\phi\rangle\,ds
\]
satisfies
\(\sup_{t\in[0,T]}|\widetilde M^{(n),\phi}_t-M^{(n),\phi}_t|\le T\varepsilon_n\).
Moreover, since \(|\langle\mu^{(n),k}_t,\phi\rangle|\le\|\phi\|_\infty\), the
bounds \eqref{eq:gamma1-bound-tightness} and \eqref{eq:gamma2-bound-tightness}
give, deterministically and uniformly in \(n\)
\begin{equation}
\label{eq:deterministic-martingale-bounds}
\sup_{t\in[0,T]}\big|M^{(n),\phi}_t\big|
\le 2\|\phi\|_\infty+C_\phi T ,
\qquad
\langle M^{(n),\phi}\rangle_T\le C_\phi T .
\end{equation}
Finally, by Corollary~\ref{cor:qv-interaction-vanishes} we have the following
convergence in \(L^1\)
\begin{equation}
\label{eq:bracket-error-identification}
\sup_{t\in[0,T]}
\left|
\langle M^{(n),\phi}\rangle_t
-
\frac{\alpha}{k}
\int_0^t
\langle \mu^{(n),k}_s,\nabla\phi\cdot a\nabla\phi\rangle\,ds
\right|
\longrightarrow 0 .
\end{equation}

Let \(0\le s<t\le T\), and let
\(G\colon D([0,s],\mathcal M_1(\mathbb R^d))\to \mathbb R\) be bounded and
continuous. We write \(G^{(n)}:=G\big((\mu^{(n),k}_u)_{0\le u\le s}\big)\). Since
\(M^{(n),\phi}\) and \((M^{(n),\phi})^2-\langle M^{(n),\phi}\rangle\) are
martingales, we have
\begin{align*}
\mathbb E\left[
\left(M^{(n),\phi}_t-M^{(n),\phi}_s\right)G^{(n)}
\right]&=0,\\
\mathbb E\left[
\left(
(M^{(n),\phi}_t)^2-(M^{(n),\phi}_s)^2
-\langle M^{(n),\phi}\rangle_t+\langle M^{(n),\phi}\rangle_s
\right)G^{(n)}
\right]&=0 .
\end{align*}
By the bound on \(\widetilde M^{(n),\phi}-M^{(n),\phi}\) and by
\eqref{eq:deterministic-martingale-bounds} and
\eqref{eq:bracket-error-identification}, both identities remain valid, up to
an error which vanishes as \(n\to\infty\), if \(M^{(n),\phi}\) is replaced by
\(\widetilde M^{(n),\phi}\) and \(\langle M^{(n),\phi}\rangle_t\) by
\(\frac{\alpha}{k}\int_0^t\langle \mu^{(n),k}_u,\nabla\phi\cdot
a\nabla\phi\rangle\,du\). The resulting functionals of \(\mu^{(n),k}\) are
bounded uniformly in \(n\), since \(\phi\), \(\mathcal A\phi\) and
\(\nabla\phi\cdot a\nabla\phi\) are bounded and all the measures have total
mass one. Moreover, since \(\mu\) has continuous paths by
Theorem~\ref{thm:tightness-empirical-measures}, they are continuous at \(\mu\)
with respect to the Skorokhod topology. Hence, using the convergence
\(\mu^{(n),k}\Rightarrow\mu\), we obtain
\begin{align*}
\mathbb E\left[
\left(M^\phi_t-M^\phi_s\right)
G\big((\mu_u)_{0\le u\le s}\big)
\right]&=0,\\
\mathbb E\left[
\Big(
(M^\phi_t)^2-(M^\phi_s)^2
-\frac{\alpha}{k}\int_s^t
\langle \mu_u,\nabla\phi\cdot a\nabla\phi\rangle\,du
\Big)
G\big((\mu_u)_{0\le u\le s}\big)
\right]&=0,
\end{align*}
where
\(M^\phi_t:=\langle \mu_t,\phi\rangle-\langle \mu_0,\phi\rangle
-\int_0^t\langle \mu_u,\mathcal A\phi\rangle\,du\) is the process of
\eqref{eq:DK-martingale}.

By a monotone class argument both identities extend from \(G\) to all bounded
random variables measurable with respect to the past of \(\mu\) up to time
\(s\). Hence \(M^\phi\) and
\((M^\phi_t)^2-\frac{\alpha}{k}\int_0^t\langle \mu_u,\nabla\phi\cdot
a\nabla\phi\rangle\,du\) are martingales. Since \(\mu\) has continuous paths,
\(M^\phi\) is continuous, and therefore its quadratic variation is given
by \eqref{eq:DK-bracket}.

Finally, since the initial measures are deterministic and
\(\mu^{(n),k}_0\to\mu^k_0\), the initial coordinate of any subsequential limit
satisfies
\[
\mu_0=\mu^k_0
\qquad\text{almost surely}.
\]
Therefore \(\mu\) solves the Dean-Kawasaki martingale problem of
Definition~\ref{def:DK-mp}.
\end{proof}

\begin{proof}[Proof of Theorem~\ref{thm:empirical-to-dk}]
The family \((\mu^{(n),k})_{n\ge1}\) is tight with continuous limit paths by
Theorem~\ref{thm:tightness-empirical-measures}, and by
Proposition~\ref{prop:identification-limit-points} every subsequential limit
solves the martingale problem of Definition~\ref{def:DK-mp} started at
\(\mu^k_0\).  That problem is well posed for atomic initial data by
Theorem~\ref{thm:DK-uniqueness-rigidity}, so all subsequential limits share the
law of \(\mu^k_t=\frac1k\sum_{i=1}^k\delta_{Y_i(t)}\), with \(Y_1,\dots,Y_k\)
independent diffusions with generator \(\frac{\alpha}{2}a:D^2\), and the whole
sequence converges.
\end{proof}

\subsection*{Acknowledgment}

JZ: Research was sponsored by the Army Research Office and was accomplished under Grant Number W91INF-23-1-0230. The views and conclusions contained in this document are those of the authors and should not be interpreted as representing the official policies, either expressed or implied, of the Army Research Office or the U.S.\ Government. The U.S.\ Government is authorized to reproduce and distribute reprints for
Government purposes notwithstanding any copyright notation herein.


\end{document}